\documentclass{article}

\usepackage{arxiv}

\usepackage[utf8]{inputenc} % allow utf-8 input
\usepackage[T1]{fontenc}    % use 8-bit T1 fonts
\usepackage{hyperref}       % hyperlinks
\usepackage{url}            % simple URL typesetting
\usepackage{booktabs}       % professional-quality tables
\usepackage{amsfonts}       % blackboard math symbols
\usepackage{nicefrac}       % compact symbols for 1/2, etc.
\usepackage{microtype}      % microtypography
\usepackage{graphicx}
\usepackage[square, comma, sort&compress, numbers]{natbib}
\usepackage{amsmath,amsthm,xcolor}

\newtheorem{assumption}{Assumption}
\newtheorem{remark}{Remark}[section]
\newtheorem{lemma}{Lemma}[section]
\newtheorem{theorem}{Theorem}[section]

\title{Revisiting the central limit theorems for the SGD-type methods}

\author{Tiejun Li,\thanks{tieli@pku.edu.cn}\quad  Tiannan Xiao,\thanks{alxeusxiao@pku.edu.cn}\quad and  Guoguo Yang\thanks{ygj512@hotmail.com}\\
LMAM and School of Mathematical Sciences,\\
Peking University, Beijing 100871, China
}

\renewcommand{\shorttitle}{Revisiting the CLT for the SGD-type methods}

\begin{document}
\maketitle

\begin{abstract}
We revisited the central limit theorem (CLT) for stochastic gradient descent (SGD) type methods, including the vanilla SGD,  momentum SGD and Nesterov accelerated SGD methods with constant or vanishing damping parameters.  By taking advantage of Lyapunov function technique and $L^p$ bound estimates, we established the CLT under more general conditions on learning rates for broader classes of  SGD methods compared with previous results. The CLT for the time average was also investigated, and we found that it held in the linear case, while  it was not generally true in nonlinear situation. Numerical tests were also carried out to verify our theoretical analysis.
\end{abstract}

% keywords can be removed

\section{Introduction}
We consider the problem
\begin{equation}
\begin{split}
\min_{x\in \mathbb{R}^d} f(x) := \frac{1}{S}\sum_{i=1}^S f_i(x),
\end{split}
\end{equation}
where $f, f_i:\mathbb{R}^d \rightarrow \mathbb{R}$ are continuously differentiable functions.  When $S$, usually the number of samples in machine learning, is very large, one can utilize the well-known stochastic gradient descent (SGD) iteration
\begin{equation}\label{eq:SGD1}
\begin{split}
x_{k} & = x_{k-1} - \alpha_k \nabla f_{s_k}(x_{k-1}) \\
& = x_{k-1} - \alpha_k \nabla f(x_{k-1}) + \alpha_k \xi_k
\end{split}
\end{equation}
to approximately locate the minima of $f$ under suitable convexity condition, where $\alpha_k$ is the learning rate, $s_k$ is a random variable uniformly sampled from $\{1,2,\ldots,S\}$. In \eqref{eq:SGD1}, the noise term $\xi_k=\nabla f(x_{k-1})-\nabla_{s_k} f(x_{k-1})$  satisfies the centering condition $\mathbb{E} [\xi_k | x_{k-1}] = 0$.

There has been a lot of extensions and theoretical analysis on SGD-type methods in the literature. To be precise, we denote a vanilla SGD (vSGD) with a general Markovian iteration form
\begin{equation}
x_{k} = x_{k-1} - \alpha_k F(x_{k-1},\xi_k),  \label{eq:vSGD}
\end{equation}
where $\mathbb{E}[ F(x_{k-1}, \xi_k)|x_{k-1}]=\nabla f(x_{k-1})$. Its different variants, such as the  momentum SGD (mSGD) (SGD version of Polyak's Heavy-ball method \citep{1964Some})
\begin{equation}\label{eq:mSGD}
\begin{split}
v_k & = \beta_k v_{k-1} - \alpha_k F(x_{k-1},\xi_k),\\
x_{k} & = x_{k-1} + v_k,
\end{split}
\end{equation}
or the Nesterov accelerated form  (NaSGD) \citep{1983A}
\begin{equation}\label{eq:NaSGD}
\begin{split}
x_{k} & = y_{k-1} - \alpha_k F(y_{k-1},\xi_k), \\
y_k & = x_k + \beta_k (x_k -  x_{k-1}),
\end{split}
\end{equation}
are also widely used to speed up the convergence, where $\beta_k\in [0,1)$ in both expressions.   For vSGD \eqref{eq:vSGD}, a well-known result is that the almost sure convergence of $\{x_k\}$ towards a critical point of $f$  holds when the learning rates $\{\alpha_k\}$ fulfill the typical condition $\sum \alpha_k = \infty$ and $\sum \alpha_k^2 < \infty$  even for nonconvex problems \citep{Bertsekas00}. The first non-asymptotic convergence rate was obtained in \cite{SGDZO} with the form $\min_{k \le n}\mathbb{E} \| \nabla f(x_k)\|^2 = O(1/\sqrt{n})$, which is un-improvable in general even when $f$ is (non-strongly) convex. This rate can be improved to be $O(1/n)$  when the variance reduction technique was utilized \citep{NCSGDP1}. However, for strongly convex $f$, the optimal convergence rate can be shown to be $\mathbb{E}\| \nabla f(x_n)\|^2 = O(1/n)$  (or similar estimates with respect to $f(x_n)-f(x^*)$ or $\|x_k-x^*\|^2$, where $x^*$ is the unique minima) when the learning rate $\alpha_k=O(1/k)$
\citep{BLN,SGDASGD,JMLR:v20:18-759}. By contrast, the analysis for the mSGD and NaSGD are relatively few. The asymptotic convergence of a general class of SGD with momentum was proven in \cite{MGDC}, and the convergence rate of the time averaged  $\mathbb{E} \| \nabla f(x_k)\|^2$ of mSGD was studied in \cite{MSGD2}. The stationary convergence bound of this time averaging was also established in \cite{SGDB} when the learning rates are  constant. For Nesterov acceleration gradient, it is well-known that  the convergence rate of $f(x_n)$ towards $f(x^*)$ is $O(n^{-2})$ in deterministic case when $f$ is convex \citep{NAG1}. And the stationary convergence bound of NaSGD was  studied in \cite{NASGDB} when the learning rates are constant.

Besides the aforementioned convergence analysis, it is a natural question to study the law of fluctuations  after knowing the convergence of $\{x_k\}$, which is expected to be related to the central limit theorem (CLT) in classical probability theory. However, this topic is less well studied in the literature.  We would like to mention that the CLT was investigated in \cite{CLTO} for the SGD-type methods with/without momentum. And it was obtained that  the convergence in distribution
\begin{equation}\label{eq:CLT1}
\begin{split}
\frac{x_n-x^*}{\sqrt{\alpha_n}} \Rightarrow N(0,V)
\end{split}
\end{equation}
under the condition that $f$ is strongly convex and $\sum_{k=1}^{\infty} \alpha_k^{3/2} < +\infty$, where $V$ is the asymptotic covariance matrix. { In \cite{borkar2009stochastic}, the condition had relaxed to $\sum_{k=1}^{\infty} \alpha_k^{2} < +\infty$.} The CLT for the variance reduced SGD (SVRG) can be enhanced to $(x_n-x^*)/a^n \Rightarrow N(0,V)$ for some $a\in (0,1)$ \citep{lei2020VR}. It was also shown that $(\hat{x}_n-x^*)/\sqrt{n} \Rightarrow N(0,V)$ where $\hat{x}_n := \sum_{k=1}^n x_k/n$ is the empirical average of $x_k$ \citep{ASGD}. And a continuous time version of CLT was given in \cite{sirignano2019} when $f$ is strongly convex.

The purpose of this paper is to give a systematic restudy on the CLT for the SGD-type methods. We will consider the CLT for the vSGD \eqref{eq:vSGD}, NaSGD \eqref{eq:NaSGD}, and the mSGD with the form
\begin{equation}\label{eq:mSGD2}
\begin{split}
v_k & = (1 - \beta_k) v_{k-1} - \alpha_k F(x_{k-1},\xi_k), \\
x_{k} & = x_{k-1} + \alpha_k v_k,
\end{split}
\end{equation}
where $\beta_k=\mu_k\alpha_k$ and $\mu_k>0$ is a damping parameter. We take the form \eqref{eq:mSGD2} instead of \eqref{eq:mSGD} due to its better  connection with the continuous time limit.

The main contributions of this paper are in three folds.

\begin{enumerate}

\item \textit{CLT for vSGD.} We re-establish the CLT for vSGD by borrowing the idea in \cite {sirignano2019}. But our different setup on the assumptions on noise and learning rates complicates the overall analysis, which is different from that taken in \cite{CLTO} and \cite{sirignano2019}.  With our approach, we can get the same CLT form \eqref{eq:CLT1} as in \cite{CLTO} but only require the learning rates satisfy the condition \eqref{eq:AssumAlpha} and the so-called slow condition (see Assumption \ref{Assum:d-slower}). This  relieves the conditions on $\{\alpha_k\}$ in \cite{CLTO}, which requires $\sum_k \alpha_k ^{2} < \infty$ in linear case and   $\sum_k \alpha_k ^{3/2} < \infty$ in nonlinear case (i.e., $\nabla f(x)$ is linear or nonlinear on $x$).

\item \textit{CLT for mSGD and NaSGD.} This part can be classified into two cases: the case with constant damping $\mu_k\equiv\tilde{\mu}$, or the case with vanishing damping $\mu_k\rightarrow 0$.  In the continuum limit, this corresponds to the SDEs
\begin{equation}
\ddot{x}+\mu(t) \dot{x}+\nabla f(x)+\xi=0
\end{equation}
with constant or vanishing damping $\mu(t)$ \citep{NAG1}. Taking advantage of the Lyapunov function technique, we
can show the CLT holds for both cases. The CLT for the vanishing damping case has never been studied before, and the constant damping case was only implicitly considered for the mSGD  in the linear case in \cite{CLTO},  but with stronger requirement $\sum_k \alpha_k ^{2} < \infty$. {Besides,  a two-time-scale stochastic approximation was given in \cite{borkar2009stochastic},
\begin{equation*}\label{eq:twoscale}
\begin{split}
w_k & = w_{k-1} - \alpha_k (A_{11} w_{k-1} + A_{12}u_{k-1} + \xi_k), \\
u_{k} & = u_{k-1} - \beta_k (A_{21} w_{k-1} + A_{22}u_{k-1} + \eta_k),
\end{split}
\end{equation*}
where $\alpha_k = o(\beta_k)$, which is quite different from the mSGD shown in Section \ref{sec:msgd}.}

\item \textit{CLT for the time average of $x_k$.} We discussed the CLT for the time average:
\begin{equation} \label{eq:aSGD}
\bar{x}_n = \frac{\sum_{k=1}^n \alpha_k x_{k-1}}{\sum_{k=1}^n \alpha_k},
\end{equation}
which is the discrete analog of the continuous time average $\int_0^T x(t)dt /T$, where $T$ is the summation of step size. Interestingly, we found that  the CLT of the form $T_n/\sqrt{S_n} \cdot(\bar{x}_n - x^*) \Rightarrow N(0,V)$ holds in the linear case, while  it is not true in general nonlinear case (see details in Sec. \ref{sec:average}).
\end{enumerate}

The rest of this paper is organized as follows. We will prove the CLT for the vSGD, mSGD and NaSGD, and the time average $\bar{x}_n$ in Sections \ref{sec:sgd}, \ref{sec:msgd} and \ref{sec:average},  respectively. Then we present numerical examples to verify the obtained CLT results in Section \ref{sec:numer}. Finally, we make the conclusion. Some proof details are left in the Appendix.

\section{Assumptions and Lemmas}
Let us first present some preliminary mathematical setup utilized in this paper. We will always assume the following necessary condition on $\{\alpha_k\}$ to ensure the convergence of $\{x_k\}$ to $x^*$ when $f$ is strongly convex.
\begin{assumption} Assumptions on the learning rates.
\begin{enumerate}
\item[(A1)] ($\mathbf{Divergence ~condition}$)\label{Assum1:alpha}
The learning rates  $\{\alpha_k\}$ satisfy
\begin{equation}\label{eq:AssumAlpha}
\lim_{k\rightarrow \infty}\alpha_k=0,\quad\sum_{k=1}^\infty \alpha_k=\infty.
\end{equation}

\item[(A2)] ($h_0$-$\mathbf{slow ~condition} )$\label{Assum:d-slower}
The learning rates $\{\alpha_k\}$ are said to satisfy the $h_0$-slow condition with $h_0>0$ if $\{\alpha_k\}$ fulfill Assumption (A1) and
\begin{equation}\label{eq:d-slower}
 \exists K_s>0, \text{ such that }\quad \frac{\alpha_n}{\alpha_m} \ge K_s \prod_{k=m + 1}^n (1 - h_0 \alpha_k) \quad \text{ for any }\ n > m \text{ large enough}.
\end{equation}

\item[(A3)]  ($\mathbf{Sufficient ~decrease ~condition}$)\label{Assum:decr}
 The learning rates $\{\alpha_k\}$ are said to be sufficiently decreasing  if $\{\alpha_k\}$ fulfill  Assumption (A1) and
\begin{equation}
 \exists d_0\ge 0, \text{ such that }\quad \frac{\alpha_k - \alpha_{k+1}}{\alpha_{k}} = d_0\alpha_k + o(\alpha_{k})\quad \text{for large enough } k.
\end{equation}
Indeed we have $d_0=\lim_{k\rightarrow\infty}(\alpha_k-\alpha_{k+1})/\alpha_k^2\ge 0$ if $\{\alpha_k\}$ satisfy sufficient decrease condition.
\end{enumerate}
\end{assumption}

Note that if $\{\alpha_k\}$ satisfy the $h_0$-slow condition, then they also satisfy $h$-slow condition for any $h > h_0$.
The proposed $h_0$-slow condition is weak enough to cover some commonly used learning  schedules with slower decreasing rate than $O(1/k)$, while strong enough to build up our CLT estimates.

It is straightforward to observe that if $\{\alpha_k\}$ are sufficiently decreasing, they will fulfill the $h_0$-slow condition for any $h_0>d_0$ since
\[
\frac{\alpha_n}{\alpha_m} = \prod_{k = m}^{n-1}\left(1 - d_0\alpha_k - o(\alpha_k)\right) \ge \frac12\prod_{k = m+1}^n(1 - h_0\alpha_k)\quad \text{for large enough }m.
\]

Typical examples which satisfy $h_0$-slow condition and sufficient decrease condition include:
\begin{enumerate}
\item[1)] $\alpha_k = K/ k^a, K>0, a\in (0,1)$, then $d_0 =0$.
\item[2)] $\alpha_k = K/k, K>0$, then $d_0 =K^{-1}$.
\item[3)] $ \beta_k=k\alpha_k$  increases to infinity monotonically, e.g. $\alpha_k = Ck^{-a} \ln k, a\in (0,1]$, then
$$d_0=\lim _{k \rightarrow \infty} \alpha_k^{-1}\left(1-\frac{k \beta_{k+1}}{(k+1) \beta_k}\right) \leq \lim _{k \rightarrow \infty} \frac{1}{(k+1) \alpha_k}=0.$$
\item[4)] If $\{\alpha_k\}$ satisfy the $h_0$-slow condition and $\{\eta_k\}$ is a bounded positive sequence, i.e., $0 < \inf_k {\eta_k} \le \sup_k {\eta_k} < +\infty$,  then $\{\alpha_k\eta_k\}$ satisfy the $h_0/\inf_k {\eta_k}$-slow condition.
\end{enumerate}

Define the filtration
\begin{equation}
\mathcal{F}_k = \sigma(x_0,\xi_1,\xi_2,\cdots,\xi_k),
\end{equation}
i.e., the $\sigma$-algebra generated by $(x_0,\xi_1,\ldots, \xi_k)$.  Denote by $X_n \overset{p}{\rightarrow} X$ the convergence in probability, i.e., $\lim_{n\rightarrow \infty}\mathbb{P}(\|X_n-X\|\ge \varepsilon)=0$ for any $\varepsilon>0$, where $\|\cdot\|$ could be any norm for $X_n,X\in \mathbb{R}^d$ or $\mathbb{R}^{d\times d}$ because of the norm equivalence theorem in finite dimensions.

\begin{assumption}Statistics for $\{\xi_k\}$ and initial value.
\begin{enumerate}
\item[(A4)]
\label{Assum:StatXi}
 We assume the following conditional mean and covariance conditions for the noise term $\{\xi_k\}$ in SGD-type methods
\begin{equation}\label{eq:StatXi}
\mathbb{E}[\xi_k|\mathcal{F}_{k-1}] = 0,\quad \mathbb{E} [\xi_k \xi_k^T|\mathcal{F}_{k-1} ] \overset{p}{\rightarrow}\Sigma \succ 0,
\end{equation}
where the notation `$A\succ B$' means $A-B$ is a symmetric positive definite (SPD) matrix.

\item[(A5)] For the initial value $x_0$, we assume that there exists large enough $p\ge 2$, such that  $\mathbb{E}\| x_0 \|^p < \infty$ and $\mathbb{E}|f(x_0)| < \infty$.

\item[(A6)]  Assumption (A4), (A5) hold, and  there exist $\delta_0,M,K_\xi>0$ such that
 \begin{equation}\label{eq:XiBound}
\mathbb{E}\big[\| \xi_k \|^{2+\delta_0}|\mathcal{F}_{k-1}\big] \le  M + K_{\xi} \| \nabla f(x_{k-1}) \|^{2+\delta_0}.
\end{equation}
\end{enumerate}
\end{assumption}

\begin{remark}
 Assumption (A4) on $\left\{\xi_k\right\}$ holds naturally for the incremental SGD form (\ref{eq:SGD1}) in practice. According to \cite{Bertsekas00}, one can establish the convergence of SGD by supplementing the condition like $\left\|\nabla f_i(x)\right\| \leq C(1+\|\nabla f(x)\|)$ for any $i$ and $x$. With strong convexity of $f$, we have $x_k \rightarrow x^*$, thus
$$
\mathbb{E}\left[\xi_k \xi_k^T \mid \mathcal{F}_{k-1}\right]=\operatorname{Cov}_{s_k}\left(\nabla f_{s_k}\left(x_{k-1}\right)\right) \rightarrow \operatorname{Cov}_{s_k}\left(\nabla f_{s_k}\left(x^*\right)\right)=: \Sigma,
$$
where $\operatorname{Cov}_{s_k}$ means the covariance respect to $s_k \sim \operatorname{Uniform}\{1,2, \ldots, S\}$. See also Remark  \ref{remark:2} for related discussions.
\end{remark}

\begin{remark}\label{remark:2} Condition (A6) is an extension of the assumption on $\{\xi_k\}$ to establish the convergence of SGD \cite{Bertsekas00}, which requires slightly stronger moment bounds for CLT. We note that it holds for the incremental SGD form \eqref{eq:SGD1} once we assume the condition $\|\nabla f_i(x)\|\le C(1+\|\nabla f(x)\|)$ for any $i$ and $x$, which is the same as that taken in \cite{Bertsekas00}. So it is not an over-stringent condition.
Condition (A9) is automatically true if $f\in C^3$ in a neighborhood of the origin.
\end{remark}

For the objective function $f$, we usually assume the following conditions.
\begin{assumption} Convexity and regularity on  function $f$.
\begin{enumerate}
\item[(A7)]
($L$-$\mathbf{smooth}$)\label{Assum:L-smooth}
The function $f: \mathbb{R}^d\rightarrow \mathbb{R}$ is said to be $L$-smooth if $f$ is differentiable and  there exists $L > 0$, such that
\begin{equation} \label{eq:L-smooth}
\forall x,y\in \mathbb{R}^d, \quad \|\nabla f(x) - \nabla f(y) \| \le L \| x - y \|.
\end{equation}

\item[(A8)] ($\mu$-$\mathbf{strongly ~convex}$)\label{Assum:StrongConv}
The function $f: \mathbb{R}^d\rightarrow \mathbb{R}$ is called $\mu$-strongly convex if there exists $\mu> 0$, such that
\begin{equation}\label{eq:muStrong}
\forall x,y\in \mathbb{R}^d, \quad f(x) \ge f(y) + \nabla f(y)^T(x-y) + \frac{\mu}{2}\| x - y \|^2.
\end{equation}

\item[(A9)]  Assumptions (A7), (A8) hold, and there exist  $p_0, r_0, K_d > 0$, such that
\begin{equation}
\| \nabla f(x) - \nabla^2 f(x^*)(x-x^*) \| \le K_d\| x-x^* \| ^ {1+p_0}\qquad  \text{ for any } \| x-x^*\| < r_0.
\end{equation}
\end{enumerate}
\end{assumption}

If $f$ is $L$-smooth, we can easily obtain
\begin{equation} \label{eq:L-smooth2}
\forall x,y\in \mathbb{R}^d, \quad f(x) \le f(y) + \nabla f(y)^T(x-y) + \frac{L}{2}\| x - y \|^2.
\end{equation}
And if $f$ is convex and twice differentiable further, we have $ 0\preceq\nabla^2 f(x)\preceq L\cdot I$ for any $x$, where  $A\preceq B$ means that the matrix $B-A$ is symmetric positive semidefinite (SPSD).

The $\mu$-strong convexity is equivalent to $(\nabla f(x) - \nabla f(y))^T(x - y)  \ge \mu \| x - y \|^2$, or $\nabla^2 f(x)\succeq \mu I$ when $f$ is twice differentiable, and has the implication  $\|\nabla f(x) - \nabla f(y) \| \ge \mu \| x - y \|$ \citep{Nesterov2004}.

The following martingale CLT is fundamental for establishing our CLT for SGD-type methods.
\begin{lemma}[Martingale difference CLT]\label{lemma:Mart-CLT}  Assume $\{Y_{nk}\}$ is a zero-mean, square-integrable martingale difference triangular array, i.e., $Y_{nk}\in \mathbb{R}^d$ is $\mathcal{F}_k$-adapted and $\mathbb{E} [Y_{nk} | \mathcal{F}_{k-1}] = 0$ for $n\in \mathbb{N}$, $k= 1,2,\dots,n$. If the multivariate conditional Lindeberg and variance conditions
\begin{equation}\label{eq:MartCLT}
\forall \varepsilon > 0, \ \sum_{k=1}^n \mathbb{E} \left[\|Y_{nk}\|^21_{\{\| Y_{nk} \| \ge \varepsilon \}}| \mathcal{F}_{k-1}\right] \overset{p}{\rightarrow}  0,\quad\text{and}\quad \sum_{k=1}^n\mathrm{Cov}  (Y_{nk}|\mathcal{F}_{k-1}) \overset{p}{\rightarrow} \Sigma
\end{equation}
hold, where $\Sigma$ is a fixed SPSD matrix, then $\sum_{k=1}^n Y_{nk} \Rightarrow N(0,\Sigma)$.
\end{lemma}

The one-dimensional version of the above lemma is a classical result in \cite{1980Martingale} (Corollary 3.1 in p. 58). Its multivariate extension is also true by considering $y^a_{nk} := a^T Y_{nk}$ for any $a\in \mathbb{R}^d$, and noting the conditions
\begin{equation*}
\begin{split}
\sum_{k=1}^n\mathrm{Cov} (y^a_{nk}|\mathcal{F}_{k-1}) & \overset{p}{\rightarrow}  a^T \Sigma a,
\end{split}
\end{equation*}
and for any induced matrix norm $\|\cdot\|$,
\begin{equation*}
\begin{split}
\sum_{k=1}^n\mathbb{E} \left[(y^a_{nk})^2 1_{\{ |y^a_{nk}| \ge \varepsilon\}} | \mathcal{F}_{k-1} \right]& \le \| a \|^2 \sum_{k=1}^n\mathbb{E} \left[\|Y_{nk} \|^2 1_{\{\| Y_{nk} \| \ge \varepsilon/\|a\|\}} | \mathcal{F}_{k-1} \right]\overset{p}{\rightarrow}  0.
\end{split}
\end{equation*}
We thus have $\sum_{k=1}^n y^a_{nk} \Rightarrow N(0,a^T\Sigma a)$ by
 one-dimensional version of the above Lemma \ref{lemma:Mart-CLT}, and then get $\sum_{k=1}^n Y_{nk} \Rightarrow N(0,\Sigma)$ by characteristic function approach to weak convergence and the arbitrariness of $a$ \citep{Book:Varadhan}.

The following lemma is important for us to establish the CLT for SGD with momentum.
\begin{lemma}[Lyapunov theorem]\label{Lyapunov:lemma2} Suppose that  $B$ is a Hurwitz matrix,  i.e., the real parts of all eigenvalues of $B$ are negative. Then for any SPSD matrix $Q$, there exists  a unique SPSD matrix $P$ such that
\begin{equation}
B^TP + PB = -Q\qquad \text{(Lyapunov equation)},
\end{equation}
and if $Q$ is SPD, then $P$ is also SPD.
\end{lemma}
\proof[Proof of Lemma \ref{Lyapunov:lemma2}]
Let $P = \int_0^{\infty} e^{tB^T}Q e^{tB} \mathrm{d}t$, which is well-defined since $\mathrm{max}\{\Re (\lambda)| \lambda \in \sigma(B)\} < 0$ for the spectrum $\sigma(B)$, where $\Re(z)$ means the real part of $z\in \mathbb{C}$.  We have
\[B^TP + PB = \int_0^{\infty}  \mathrm{d}(e^{tB^T}Q e^{tB}) = -Q.\]
If $V$ is the solution for $B^T V + VB = 0$, then we define $W(t) = e^{tB^T}V e^{tB}$. We have $W(0) = V$, $W(\infty) = 0$, and $W'(t) = e^{tB^T}(B^T V + VB)e^{tB} = 0$. So $V = W(0) = W(\infty) = 0$.
The positive definiteness of $P$ is straightforward from its definition.
\endproof

\noindent\textbf{Notational convention}. Throughout this paper, we will use $C$ as generic $O(1)$ positive constants in different bound estimates, whose value may change in different places. The following convention
\[\sum_{k=n+1}^n y_k:=0,\quad \prod_{k=n+1}^n y_k:=1\]
is also adopted for any sequence $\{y_k\}$. We also use the notation $a_n = \Theta(b_n)$ which means that there exists $ C_1, C_2 > 0$, such that $C_1 b_n\le |a_n| \le C_2 b_n$.

\section{CLT for vSGD}\label{sec:sgd}

In this section, we will re-establish the CLT for vSGD by borrowing the idea in \cite {sirignano2019}, which can relieve the conditions on $\{\alpha_k\}$ in \cite{CLTO}. However, our setup on the assumptions on noise and learning rates is more general than \cite{sirignano2019}, which makes the overall analysis more complicate.

Assume $f$ is twice differentiable, $L$-smooth (\ref{eq:L-smooth}), $\mu$-strongly convex (\ref{eq:muStrong}) and $x^*$ is the unique minima. Define
\begin{equation}
A = \nabla^2 f(x^*),\ X_k=x_k-x^*,\ R_k = A(x_{k-1} - x^*) - \nabla f(x_{k-1}).
\end{equation}
Then the vSGD iteration \eqref{eq:SGD1} can be rewritten as
\begin{equation}\label{eq:vSGDX}
X_k = (I - \alpha_k A)X_{k-1} + \alpha_k R_k + \alpha_k \xi_k.
\end{equation}

Without loss of generality, we will assume $\alpha_k \le L^{-1}$ for all $k$ since we are only concerned with the asymptotic properties of $x_k$. The $L$-smoothness of $f$ ensures  that $(I - \alpha_k A)$ is invertible. Define $B_k = \prod_{i=1}^k (I - \alpha_i A)$. Timing $B_k^{-1}$ to both sides of \eqref{eq:vSGDX} and taking summation from $k=1$ to $n$,  we get
\begin{equation}
X_n= B_n X_{0} + \sum_{k=1}^n B_n B_k^{-1} \alpha_k R_k + \sum_{k=1}^n B_n B_k^{-1} \alpha_k\xi_k.
\end{equation}

Our main result in this section is the CLT for $X_n=x_n-x^*$ under a suitable scaling.
The main theorem in this section is as below.
\begin{theorem}\label{thm:CLT-vSGD} Under the Assumptions  (A4),  (A7)--(A9) and
 Assumptions (A1), (A2)  hold with $h_0 < 2\mu$, and there exists $K_\alpha>0$, such that
\begin{equation}\label{eq:StepSize}
\alpha_n \le K_\alpha \min_{k \le n+1}\alpha_k\qquad \text{ for all } n,
\end{equation}
we have
\begin{equation}
\mathrm{Cov}(X_n)^{-1/2}(x_n - x^*) \Rightarrow N(0,I).
\end{equation}
Furthermore, if Assumption (A3) holds for  $d_0 < 2\mu$, then there exists  a unique SPD matrix $W^*$ such that
\[AW^* + W^*A^T  -d_0 W^*= \Sigma,\]
and we have
\begin{equation}
\frac{x_n - x^*}{\sqrt{\alpha_n}} \Rightarrow N(0,W^*).
\end{equation}
\end{theorem}

The proof of Theorem~\ref{thm:CLT-vSGD} relies on the following lemmas concerning the covariance and   $L^p$ bounds of $X_n$.

 \begin{lemma}\label{lm:lm3}
Consider the triplet $(A_k, \beta_k, W)$ with the condition $\beta_k \rightarrow 0$, $W = W^T$, and $A_k$ satisfies the recursive relation
\[A_k = (I - \alpha_k D + o(\alpha_k)) A_{k-1},\]
where $-D$ is a Hurwitz matrix.  Further assume  the learning rates satisfy  Assumption (A1),  then we have
\begin{equation}\label{eq:lemma3}
\sum_{k=1}^{n} A_n A_{k}^{-1}W A_k^{-T} A_n^T \alpha_k \beta_k \rightarrow 0.
\end{equation}
\end{lemma}

\begin{lemma}\label{lm:vSGD}
Assume that $f$ is $L$-smooth and $\mu$-strongly convex, and Assumption (A6) and  (\ref{eq:StepSize}) in Theorem~\ref{thm:CLT-vSGD}  hold. We have the following  upper bounds about $x_n-x^*$:
\begin{equation}\label{eq:UBound}
\mathbb{E}\| x_n - x^* \|^{2p} \le U_p \alpha_n^p, \quad \text{ for }\quad p \in \Big[1, 1+\frac{\delta_0}{2}\Big],
\end{equation}
where $\delta_0$ is the constant defined in \eqref{eq:XiBound}.
\end{lemma}
The proof of these two lemmas will be deferred to Appendix \ref{sec:A1}-\ref{sec:A2}. Indeed, there are similar $L^{2p}$ lower bounds for $x_n-x^*$. Lemma~\ref{lm:vSGD} tells us that the convergence of $x_n$ to $x^*$ is $O(\sqrt{\alpha_n})$ under suitable conditions on the learning rates. In fact, if some of these conditions  are not valid, e.g. $\alpha_n = K n^{-1}$ for $K>2\mu$, which violates the slow condition, it is possible that the upper bound is not $O(\sqrt{\alpha_n})$. But we will not pursuit this point in the current paper.

\proof{Proof of Theorem~\ref{thm:CLT-vSGD}}
We will prove the linear case at first, then the nonlinear case. In linear case, the vSGD has the form
\begin{equation}\label{eq:X_nLinear}
X_n =  B_nX_0+\sum_{k=1}^n B_nB_k^{-1} \alpha_k \xi_k,\quad  B_n = \prod_{k=1}^n (I - \alpha_k A).
\end{equation}
Since $B_nX_0\rightarrow 0$ almost surely, we will skip this term in the analysis below. The proof in linear case will be accomplished in four steps.

\vspace*{0.5em}
\textbf{Step 1.} Strong convergence of the conditional covariance of $\xi_k$.

First note that the condition \eqref{eq:XiBound} and the upper bound \eqref{eq:UBound}  yields  the estimate $\mathbb{E} \|\xi_k\|^{2+\delta_0}\le C$ for any $k$, where $C$ is a constant depending on $K_\xi,M,L,\delta_0$. Combining this uniform $L^{2+\delta_0}$-norm bound on $\xi_k$ and  the convergence-in-probability condition \eqref{eq:StatXi}, we can easily get
\begin{equation}\label{eq:CovLim}
\Sigma_k(\omega):=\mathbb{E}[\xi_k \xi_k^T| \mathcal{F}_{k-1} ]\rightarrow  \Sigma \qquad \text{in } L^{1,1}
\end{equation}
by splitting estimates on the domain  $\{\|\Sigma_k(\omega)-\Sigma\|\ge \varepsilon\}$ and $\{\|\Sigma_k(\omega)-\Sigma\|< \varepsilon\}$ separately,  where the $L^{1,1}$-norm $\|\Sigma\|_{1,1}:=\sum_{i,j=1}^d |\Sigma_{ij}|$ is the entry-wise $L^1$-norm of matrix $\Sigma$. Indeed, this could be replaced with any norm due to the norm equivalence theorem. Define $\Sigma_k=\mathbb{E}(\xi_k \xi_k^T)$, then we have $\Sigma_k \rightarrow \Sigma$ by \eqref{eq:CovLim}. Without loss of generality, we can assume there exists $\Gamma\succ 0$, and $\Sigma_k\succ \Gamma$ since only the limit behavior matters. By Lyapunov theorem, there exists $W\succ 0$ such that $AW + W A^T = \Gamma$.

\vspace*{0.5em}
\textbf{Step 2.} Covariance estimate.

According to \eqref{eq:StepSize}, we get
\begin{align}
\mathrm{Cov}(X_n) = \sum_{k=1}^n  B_nB_k^{-1} \Sigma_k B_k^{-T} B_n^T \alpha_k^2 \succeq {K_\alpha^{-1}}\alpha_n\sum_{k=1}^n  B_n B_k^{-1} \Gamma B_k^{-T} B_n^T\alpha_k.
\end{align}
For each summation term, we have
\begin{align}
\ B_n B_k^{-1} \Gamma & B_k^{-T} B_n^T\alpha_k = B_n B_k^{-1}(\alpha_k AW + \alpha_k W A^T ) B_k^{-T} B_n^T \nonumber \\
= & \ B_n B_k^{-1}(W - (I - \alpha_k A) W(I - \alpha_k A^T) + \alpha_k^2 AWA^T )B_k^{-T} B_n^T \nonumber \\
= &\ B_n B_{k}^{-1}WB_k^{-T} B_n^T - B_n B_{k-1}^{-1}WB_{k-1}^{-T} B_n^T + \alpha_k^2  B_n B_{k}^{-1}AWA^T B_k^{-T} B_n^T. \label{eq:CovEst1}
\end{align}
The summation of the last term in \eqref{eq:CovEst1} converges to $0$ by Lemma \ref{lm:lm3}. We get
\begin{align}
\mathrm{Cov}(X_n) & \succeq \frac{\alpha_n}{K_\alpha}\left[\sum_{k=1}^n (B_n B_{k}^{-1}WB_k^{-T} B_n^T - B_n B_{k-1}^{-1}WB_{k-1}^{-T} B_n^T) + o(1)\right] \nonumber\\
& \succeq \frac{\alpha_n}{K_\alpha}\Big(W - B_n W B_n^T + o(1)\Big) \sim C \alpha_nI \label{eq:3.13}
\end{align}
since $B_nWB_n^T\rightarrow 0$ almost surely, where $C$ depends on $K_\alpha, W$.

\vspace*{0.5em}
\textbf{Step 3.} Conditional covariance estimate.

To apply the martingale difference CLT, we define
\begin{equation}
y_{nk} = \alpha_k B_nB_k^{-1} \xi_k, \quad Y_{nk} = \mathrm{Cov}(X_n)^{-1/2} y_{nk}.
\end{equation}
To verify the conditional variance condition $\sum_{k=1}^n \mathrm{Cov}({Y_{nk}|\mathcal{F}_{k-1}})  \overset{p}{\rightarrow} I$, it is enough to show
\begin{equation}\label{eq:3.15}
G_n:=\frac1{\alpha_n}\left(\sum_{k=1}^n \mathrm{Cov}(y_{nk}|\mathcal{F}_{k-1}) -  \mathrm{Cov}(X_n) \right) \overset{p}{\rightarrow} 0
\end{equation}
by \eqref{eq:StepSize} and \eqref{eq:3.13}. The left hand side of \eqref{eq:3.15} has the form
\begin{equation}
G_n=\frac{1}{\alpha_n} \left( \sum_{k=1}^n B_nB_k^{-1}\Lambda_k B_k^{-T}B_n^T\alpha_k^2 \right),\quad \Lambda_k :=  \mathbb{E}[\xi_k \xi_k^T|\mathcal{F}_{k-1}]  - \Sigma_k
\end{equation}
where $\Lambda_k \rightarrow 0$ in $L^{1,1}$ according to Step 1. Define
\[A_n = \prod_{k=1}^n (I - \alpha_k A)(1- h_0 \alpha_k)^{-1/2} = \prod_{k=1}^n (I - \alpha_k (A - h_0 / 2) + \beta_k),\quad \beta_k = o(\alpha_k).\]
We have
\begin{align}\label{eq:3.17}
\| G_n  \|_2 \le &   \left\| \sum_{k=1}^n B_nB_k^{-1} B_k^{-T}B_n^T\frac{\alpha_k^2}{\alpha_n} \| \Lambda_k \|_2 \right\|_2   \le K_s ^{-1} \left\| \sum_{k=1}^n A_nA_k^{-1}A_k^{-T}A_n^T {\alpha_k} \| \Lambda_k \|_2  \right\|_2.
\end{align}
The first inequality in \eqref{eq:3.17} is from the fact that $\Lambda_k\preceq \|\Lambda_k\|_2 I$ and $\|P\|_2\le \|Q\|_2$ for any symmetric $P,Q$ which satisfy $Q\succeq 0, Q \succeq P \succeq -Q$.

To prove that  $\mathbb{E}\| G_n  \|_2\rightarrow 0$, we only need to show that
\[F_n:=\sum_{k=1}^n A_nA_k^{-1}A_k^{-T}A_n^T {\alpha_k} \mathbb{E}\| \Lambda_k \|_2\rightarrow 0\]
since each summand in $F_n$ is a SPSD matrix. The convergence $F_n\rightarrow 0$ can be established by Lemma~\ref{lm:lm3} by choosing the triplet $(A_k, \mathbb{E}\|\Lambda_k\|_2, I)$.

The variance condition is verified.

\vspace*{0.5em}
\textbf{Step 4.} Conditional Lindeberg condition.

To establish the  conditional Lindeberg condition \eqref{eq:MartCLT}, first we note that
\begin{equation*}
\begin{split}
\sum_{k=1}^n \mathbb{E} \big[\| Y_{nk} \|^2 1_{\{\| Y_{nk} \| \ge \varepsilon\}}|\mathcal{F}_{k-1}\big] &  \le \varepsilon^{-\delta} \sum_{k=1}^n  \mathbb{E}  [\| Y_{nk} \|^{2+\delta}|\mathcal{F}_{k-1}]
\end{split}
\end{equation*}
fo any $\delta\ge 0$. So in order to prove the convergence in probability of the above term to 0, we only need to show   $\sum_{k=1}^n  \mathbb{E}  \| Y_{nk} \|^{2+\delta} \rightarrow 0$ for some $\delta>0$.

With the estimate $\mathrm{Cov}(X_n) \succeq C\alpha_n I$, we have for $\delta<\delta_0$
\[\mathbb{E} \| Y_{nk} \|^{2+\delta} \le C \frac{\| B_n B_k^{-1} \alpha_k \|^{2+\delta}}{\alpha_n^{1+\delta/2}},\]
where $C$ depends on  $L$, $\sup_k \mathbb{E}\| \xi_k \|^{2+\delta}<\infty$ by \eqref{eq:XiBound}, and the upper bounds in Lemma~\ref{lm:vSGD}.
With  (\ref{eq:StepSize}) and $\mu$-strong convexity $A\succeq \mu I$, we obtain
\begin{equation*}
\begin{split}
\sum_{k=1}^n \| B_n B_k^{-1} \alpha_k \|^{2+\delta}  & \le (K_s^{-1}\alpha_n)^{1+\delta}\sum_{k=1}^n \left( \prod_{j = k+1}^n (1 - \mu \alpha_j)^{2+\delta} (1 - h_0 \alpha_j)^{-1- \delta}\right) \alpha_k \\
& \le (K_s^{-1}\alpha_n)^{1+\delta}\sum_{k=1}^n  \exp\left( -C_
\delta\sum_{j = k+1}^n \alpha_j \right)(1-e^{-C_\delta\alpha_k})C^{-1}_\delta e^{C_\delta\alpha_k},
\end{split}
\end{equation*}
where $C_\delta := \mu(2+\delta)-h_0(1+\delta)>0$ by choosing $0<\delta<\min\{(2\mu-h_0) / (h_0 - \mu), \delta_0\}$ if $\mu \le h_0$, or any $0<\delta<\delta_0$ otherwise. With such choice, we get
\begin{equation}\label{ineq:estimate2}
\sum_{k=1}^n \| B_n B_k^{-1} \alpha_k \|^{2+\delta} \le (K_s^{-1}\alpha_n)^{1+\delta}\left(1-   e^{ -C_
\delta\sum_{k=1}^n \alpha_k }\right)C^{-1}_\delta e^{C_\delta\alpha_k} \le
C \alpha_n^{1+\delta}.
\end{equation}
So we have
\begin{equation*}
\begin{split}
\sum_{k=1}^n \mathbb{E} \| Y_{nk} \|^{2+\delta}   \le C\frac{\alpha_n^{1+\delta}}{\alpha_n^{1 + \delta/2}} = C  \alpha_n^{\delta/2} \rightarrow 0
\end{split}
\end{equation*}
and the Lindberg condition is established.

For the nonlinear case, by Assumption (A9), for $0 < p \le \min (p_0, \delta_0+1)$, we have:
\begin{equation*}
\begin{split}
\mathbb{E}\| R_{k+1}\| &\le   K_d\mathbb{E}\Big(\| X_k \|^{1+p} 1_{\{\| X_k \| < r_0\}}\Big) + (\mu+L)\mathbb{E}\Big(\| X_k \| 1_{\{\| X_k \| \ge r_0\}}\Big) \\
& \le K_d\mathbb{E}\| X_k \|^{1+p} + (\mu+L)r_{0}^{-p}\mathbb{E}\| X_k \|^{1+p} \le C \alpha_k^{(1+p)/2}.
\end{split}
\end{equation*}
Further, let $p < 2\mu/ h_0 - 1$. With (\ref{eq:StepSize}) and $\mu$-strong convexity, we obtain
\begin{equation*}
\begin{split}
\mathbb{E}\left\| \sum_{k=1}^n B_nB_k^{-1} \alpha_k R_k\right \| & \le C\sum_{k=1}^n \prod_{j = k+1}^n (1 - \mu \alpha_j) \alpha_k\alpha_{k-1}^{(1+p)/2} \le C\sum_{k=1}^n \prod_{j = k+1}^n (1 - \mu \alpha_j) \alpha_k^{(3+p)/2} \\
& \le C\sum_{k=1}^n \prod_{j = k+1}^n (1 - \mu \alpha_j)(1 - h_0 \alpha_j)^{-(1+p)/2}\alpha_k \alpha_n^{(1+p)/2} \le C \alpha_n^{(1+p)/2},
\end{split}
\end{equation*}
where the last inequality is similar as \eqref{ineq:estimate2}. With $\mathrm{Cov}(X_n) \succeq C\alpha_n I$, we get
\[
\mathrm{Cov}(X_n)^{-1/2}\sum_{k=1}^n B_nB_k^{-1} \alpha_k R_k \overset{L^1}{\rightarrow}  0.
\]

 For brevity, we will only give the proof in linear case, i.e., $R_k=0$.
Define $V_k = \mathrm {Cov}(X_k), \Sigma_k = \mathbb{E}\xi_k \xi_k^T$, and $W_{k-1} =  V_{k-1} / \alpha_k $. It is enough to show $W_k\rightarrow W^*$ and $W^*$ satisfies
$\tilde{A}W^* + W^*\tilde{A}^T  = \Sigma$,  where $\tilde{A}=A-d_0/2I$.

By direct calculations, we get
\begin{equation}\label{eq:3.6}
V_k - V_{k-1}  = - (AV_{k-1} + V_{k-1}A^T )\alpha_k + \alpha_k^2 \Sigma_k + \alpha_k Q_k,
\end{equation}
where $Q_k = \alpha_kAV_{k-1}A^T$ is a symmetric matrix. Dividing both sides of \eqref{eq:3.6} with $\alpha_k$ and utilizing Assumption (A3), we get
\begin{equation}\label{eq:3.7}
W_k - W_{k-1} = -\alpha_{k}(AW_{k-1} + W_{k-1}A^T  - d_0 W_{k-1} -  \Sigma_k) + \alpha_k \tilde{Q}_k.
\end{equation}
where $\|\tilde{Q}_k\|=o(1)$.
Let $U_k = W_k - W^*$. We obtain
\begin{equation*}
\begin{split}
U_k = (I - \alpha_k \tilde{A}) U_{k-1} (I -\alpha_k \tilde{A}^T) + \alpha_k P_k,
\end{split}
\end{equation*}
where $P_k = \tilde{Q}_k- \alpha_k \tilde{A}U_{k-1}\tilde{A}^T+(\Sigma_k-\Sigma) = o(1)$ is symmetric.  Define $A_k=\prod_{i=1}^k (I - \alpha_i \tilde{A})$. We have
\begin{equation*}
U_n = \sum_{k=1}^n  A_nA_k^{-1} P_k A_k^{-T} A_n^T \alpha_k +  A_n U_0 A_n^T .
\end{equation*}
It is straightforward that $A_n U_0 A_n^T\rightarrow 0$ by strong convexity of $\tilde{A}$ and Assumption  (A1). Upon skipping this term, we have
\begin{equation*}
\|U_n\|_2 \le \left \| \sum_{k=1}^n  A_nA_k^{-1} A_k^{-T} A_n^T \alpha_k  \| P_k\|_2 \right\|_2
\end{equation*}
by the fact that $P_k\preceq \|P_k\|_2 I$ and $\|P\|_2\preceq \|Q\|_2$ for any symmetric $P,Q$ which satisfy $Q\succeq 0, Q \succeq P \succeq -Q$, where $\|\cdot\|_2$ is the $L^2$ norm of a matrix. So $U_n$ converges to 0 by applying Lemma~\ref{lm:lm3} to the triplet $(A_k,\|P_k\|_2, I)$. From \eqref{eq:3.7}, we also get
$AW^*+W^*A^T-d_0W^*=\Sigma$.

The analysis for general nonlinear case is similar, which only introduces an additional $o(\alpha_k^2)$ symmetric matrix term in \eqref{eq:3.6}, but does not alter the proof very much.
The proof is done.

\endproof

A special application of the above theorem is  that we have the convergence
\[
n^{a/2}(x_n - x^*) \Rightarrow N(0, \alpha_1W^*)
\]
when we take $\alpha_n = \alpha_1 n^{-a}$ for $a\in (0,1)$. This CLT with such weak conditions on the learning rates has not been obtained before in the literature.

\begin{remark}
For vSGD, we can explicitly solve the matrix $W^*$. Suppose $A$ has the eigendecomposition $A = Q^T\Lambda Q$, where $Q^TQ=I$ and $\Lambda$ is the diagonal eigenvalue matrix.  Then we have
\[
W^* = Q^T \mathcal{L}(Q\Sigma Q^T)Q,
\]
where $\mathcal{L}(\Sigma)$ is the matrix with components
\[\mathcal{L}(\Sigma)_{ij} := \frac{\Sigma_{ij}}{\lambda_i+\lambda_j-2d_0},\]
which is the solution of the Lyapunov equation $\Lambda W+W\Lambda^T-d_0W=\Sigma$ for $W$.
\end{remark}

\section{CLT  for mSGD and NaSGD}\label{sec:msgd}

In this section, we will establish the CLT for mSGD  and NaSGD in the form (\ref{eq:mSGD2}) with constant damping $\mu_k=\tilde{\mu}$ or vanishing damping $\mu_k\rightarrow 0$. It is worth noting that the theoretical results for the two cases are very few, which is a key motivation for us to study the problem with more general learning rates \eqref{eq:AssumAlpha}.

\subsection{CLT  for mSGD with constant damping $\mu_k\equiv\tilde{\mu}$}\label{sec:4.1}

Assume $f$ is twice differentiable, $L$-smooth (\ref{eq:L-smooth}), $\mu$-strongly convex (\ref{eq:muStrong}) and $x^*$ is the unique minima.
Consider the following mSGD
iteration
\begin{equation}\label{eq:mSGD31}
\begin{aligned}
&x_{k}=x_{k-1}+\alpha_{k} v_{k}, \\
&v_{k}=v_{k-1}-\mu_{k} \alpha_{k} v_{k-1}-\alpha_{k} \nabla f\left(x_{k-1}\right)+\alpha_{k} \xi_{k}.
\end{aligned}
\end{equation}
Define
\[Z_k = \left(\begin{array}{cc}
x_k - x^*\\
v_k
\end{array}\right),~~D = \left(\begin{array}{cc}
0 & -I\\
A & \tilde{\mu}I
\end{array}\right),~~E = \left(\begin{array}{cc}
A & \tilde{\mu}I\\
0 & 0
\end{array}\right).\]
Then  \eqref{eq:mSGD31} with constant damping $\mu_k=\tilde{\mu}$  can be rewritten as
\begin{equation}\label{eq:mSGD3}
\begin{split}
Z_k = ( I - \alpha_k D - \alpha_k^2 E)Z_{k-1}  + \alpha_k
\left(\begin{array}{c}
0 \\
R_{k-1} + \xi_k
\end{array}\right)+ o(\alpha_k).
\end{split}
\end{equation}
Let $B_k = \prod_{i=1}^k (I - \alpha_i D - \alpha_i^2 E)$. Without loss of generality, we can assume $B_k$ is invertible since we are only concerned with the asymptotic properties of $Z_k$. Timing $B_k^{-1}$ to both sides of \eqref{eq:mSGD3} and taking summation from $k=1$ to $n$,  we get
\begin{equation}
Z_n= B_n Z_{0}  + \sum_{k=1}^n B_n B_k^{-1}\left( \alpha_k
\left(\begin{array}{c}
0 \\
R_{k-1} + \xi_k
\end{array}\right)+o(\alpha_k)\right).
\end{equation}

 Let $\sigma(D)$ be the spectrum of $D$ and $\lambda_D:=\mathrm{min}\{\Re (\lambda)| \lambda \in \sigma(D)\}$. Our main result in this subsection is the CLT for $Z_n$ under a suitable scaling.

\begin{theorem}\label{thm:CLT-mSGD1}
 Under the Assumptions (A6), (A9) and
% \begin{itemize}
Assumptions (A1), (A2) hold with
  $h_0 < 2\lambda_D$, and there exists $K_\alpha>0$, such that
\begin{equation}\label{eq:StepSizem}
\alpha_n \le K_\alpha \min_{k \le n+1}\alpha_k\qquad \text{ for all } n,
\end{equation}
%\end{itemize}
we have
\begin{equation}
\mathrm{Cov}(Z_n)^{-1/2}Z_n  \Rightarrow N(0,I).
\end{equation}
Furthermore,  if Assumption (A3)  holds for  $d_0  < 2\lambda_D$,  there exists  a unique SPD matrix $W^*$ such that
\begin{equation}\label{eq:Lyap-mSGD}
DW^* + W^*D^T  -d_0 W^*= \tilde{\Sigma},
\end{equation}
where $\tilde{\Sigma}=\left(\begin{array}{cc}
0 & 0\\
0 & \Sigma
\end{array}\right)$. And for mSGD \eqref{eq:mSGD31}, we have
\begin{equation}
\frac{Z_n}{\sqrt{\alpha_n}} \Rightarrow N(0,W^*).
\end{equation}
\end{theorem}

\begin{remark} The constant $\lambda_D$ in Theorem~\ref{thm:CLT-mSGD1} is only for convenience of the proof. It can be replaced by a checkable constant
\[h_D:=\mathrm{min}\left\{\frac{2}{\zeta\mu},~\frac{2(1+\mu\zeta^2)}{\tilde{\mu}}\right\},\]
where $\zeta:=\mu/(2L+\tilde{\mu}^2)$. It is easy to verify that $h_D\le \lambda_D$.
\end{remark}

\proof{Proof of Theorem~\ref{thm:CLT-mSGD1}}
The proof is parallel to each step in proving Theorem \ref{thm:CLT-vSGD}.
The main difference is that we should re-establish the bounds $\mathbb{E}\| x_n - x^* \|^{2p} \le U_p \alpha_n^p$ for $p \in [1,1+\delta_0/2]$ in mSGD setup as in Lemma \ref{lm:vSGD}, which is deferred to Appendix \ref{sec:A3}. Furthermore, in Step 3 of Theorem \ref{thm:CLT-vSGD}, we have to consider the quantities with the modified form
$Y_{nk} = \mathrm{Cov}(Z_n)^{-1/2}B_n B_{k}^{-1}(0,\xi_k^T)^T$, $k = 1,2,\dots, n$, where $B_k=\prod_{i=1}^{k}( I - \alpha_i D - \alpha_i^2 E)$. The other procedures are similar.

\endproof

A special application of the above theorem is  that we have the convergence
\[
n^{a/2}Z_n \Rightarrow N(0,\alpha_1W^*)
\]
when we take $\alpha_n = \alpha_1 n^{-a}$ for $a\in (0,1)$.

\begin{remark}
For mSGD, we can explicitly solve the matrix $W^*$ when $d_0=0$. Suppose $A$ has the eigendecomposition $A = Q^T\Lambda Q$, where $Q^TQ=I$ and $\Lambda$ is the diagonal eigenvalue matrix.  Then we have
\[
W^* = \mathrm{diag}(Q^T,Q^T)\mathcal{L}(Q\Sigma Q^T)\mathrm{diag}(Q,Q),
\]
where $\mathcal{L}({\Sigma})=\left(\begin{array}{cc}
H & \tilde{B}^T\\
\tilde{B} & J
\end{array}\right)$ is the matrix with components
\[H_{ij} =  \frac{2\mu \Sigma_{ij}}{2\mu^2(\lambda_i + \lambda_j)+ (\lambda_i - \lambda_j)^2},~~
J_{ij} = \frac{\mu(\lambda_i+\lambda_j)\Sigma_{ij}}{2\mu^2(\lambda_i + \lambda_j)+ (\lambda_i - \lambda_j)^2},\]
\[
\tilde{B}_{ij} = \frac{\Sigma_{ij}(\lambda_j - \lambda_i)}{2\mu^2(\lambda_i + \lambda_j)+ (\lambda_i - \lambda_j)^2},\]
which is the solution of the Lyapunov equation $\Lambda W+W\Lambda^T=\tilde{\Sigma}$ for $W$. The simplest application is when $A=cI$ and $d_0=0$, we have
$
W^* =\left(\begin{array}{cc}
(2c\tilde{\mu})^{-1} I & 0\\
0 & (2\tilde{\mu})^{-1} I
\end{array}\right).
$
\end{remark}

\subsection{CLT  for NaSGD with constant damping $\mu_k\equiv\tilde{\mu}$}

In this subsection, we will establish the CLT for NaSGD with constant damping $\mu_k\equiv \tilde{\mu}$.  Assume $f$ is twice differentiable, $L$-smooth (\ref{eq:L-smooth}), $\mu$-strongly convex (\ref{eq:muStrong}) and $x^*$ is the unique minima. Further assume $\alpha_k/\alpha_{k-1} \rightarrow 1,$ then $\beta_k = (1 - \tilde{\mu} \alpha_k)\frac{\alpha_k}{\alpha_{k-1}} \rightarrow 1$. Consider the following NaSGD iteration
\begin{equation}\label{eq:NaSGD3}
\begin{aligned}
x_{k} & = x_{k-1} + \alpha_k v_{k}, \\
v_k & = (1 - \tilde{\mu} \alpha_k) v_{k-1}  -\alpha_k \nabla f(x_{k-1}+ \beta_k v_{k-1})  + \alpha_k \xi_k\\
& = v_{k-1} - (\tilde{\mu} + A)\alpha_k v_{k-1}  - \alpha_k A x_{k-1} +\alpha_k (R_{k-1} + \xi_k) + o(\alpha_k).
\end{aligned}
\end{equation}
where $R_k =  \nabla f(y_{k-1}) - Ay_{k-1}, y_{k-1} = x_{k-1}+ \beta_k v_{k-1}.$
Define
 \[Z_k = \left(\begin{array}{cc}
x_k - x^*\\
v_k
\end{array}\right),~~D = \left(\begin{array}{cc}
0 & -I\\
A & \tilde{\mu} I+A
\end{array}\right).
\]
Now \eqref{eq:NaSGD3} has the form
\begin{equation*} \label{eq:NaSGD31}
\begin{split}
Z_k = ( I - \alpha_k D+o(\alpha_k) )Z_{k-1}  + \alpha_k
\left(\begin{array}{c}
0 \\
R_{k-1} + \xi_k
\end{array}\right) + o(\alpha_k).
\end{split}
\end{equation*}
Let $B_k = \prod_{i=1}^k (I - \alpha_i D + o(\alpha_i))$. With similar procedure as before,  we get
\begin{equation}
Z_n= B_n Z_{0}  + \sum_{k=1}^n B_n B_k^{-1} \left(\alpha_k
\left(\begin{array}{c}
0\\
R_{k-1} + \xi_k
\end{array}\right)+o(\alpha_k)\right).
\end{equation}

Similarly define $\sigma(D)$ and $\lambda_D$ as in Theorem~\ref{thm:CLT-mSGD1} for the current matrix $D$. We have the CLT for NaSGD with constant damping.
\begin{theorem}\label{thm:CLT-NaSGD1}
Under the Assumptions (A6), (A9) and  (\ref{eq:StepSizem}) in Theorem~\ref{thm:CLT-mSGD1}, and the condition $\lim_{n\rightarrow\infty}\alpha_n/\alpha_{n-1}=1$,
 we have
\begin{equation}
\mathrm{Cov}(Z_n)^{-1/2}Z_n  \Rightarrow N(0,I).
\end{equation}
Furthermore, if Assumption  (A3) holds for  $d_0  < 2\lambda_D$, there exists  a unique SPD matrix $W^*$ satisfies \eqref{eq:Lyap-mSGD}, and for NaSGD \eqref{eq:NaSGD3}, we have
\begin{equation}
\frac{Z_n}{\sqrt{\alpha_n}} \Rightarrow N(0,W^*).
\end{equation}
\end{theorem}

\begin{remark} The constant $\lambda_D$ can be replaced by a checkable constant
\[h_D := \mathrm{min}\left\{\frac{2}{\zeta},~\frac{2(1+\mu\zeta^2)}{\tilde{\mu}}\right\},\]
where $\zeta=(\mu+\tilde{\mu})/(2L+\tilde{\mu}^2)$ since $h_D\le \lambda_D$.
\end{remark}

\proof{Proof of Theorem~\ref{thm:CLT-NaSGD1}}
The proof is similar to Theorem~\ref{thm:CLT-mSGD1} by establishing the upper bounds for $Z_n$: $\mathbb{E}\| Z_n \|^{2p} \le U_p \alpha_n^p$ for $p \in [1,1+\delta_0/2]$ in NaSGD setup, which is deferred to Appendix \ref{sec:A4}, and taking the corresponding change $Y_{nk} = \mathrm{Cov}(Z_n)^{-1/2}B_n B_{k}^{-1}(0,\xi_k^T)^T$ for $k = 1,2,\dots, n$, where $B_k=\prod_{i=1}^{k}( I - \alpha_i D + o(\alpha_i))$.

\endproof

\subsection{mSGD with vanishing damping $\mu_k \rightarrow 0$}
In this subsection, we will establish the CLT for mSGD with vanishing damping
$\mu_k \rightarrow 0$.
Assume $f$ is twice differentiable, $L$-smooth (\ref{eq:L-smooth}), $\mu$-strongly convex (\ref{eq:muStrong}) and $x^*$ is the unique minima.
Define
\[Z_k = \left(\begin{array}{cc}
x_k - x^*\\
v_k
\end{array}\right),~~\tilde{D} = \left(\begin{array}{cc}
0 & -I\\
A & 0
\end{array}\right),~~\tilde{E}= \left(\begin{array}{cc}
0 & 0\\
0 & I
\end{array}\right).\]
Then  the mSGD iteration (\ref{eq:mSGD2}) with vanishing damping $\mu_k$  has the form
\begin{equation}\label{eq:mSGD4}
\begin{split}
Z_k = \left( I - \alpha_k \tilde{D} - \alpha_k\mu_k \tilde{E}+\alpha_k^2\left(\begin{array}{cc}
A & \mu_k I\\
0 & 0
\end{array}\right)\right)Z_{k-1}  + \alpha_k
\left(\begin{array}{c}
\alpha_k (R_{k-1} + \xi_k ) \\
R_{k-1} + \xi_k
\end{array}\right).
\end{split}
\end{equation}

For the vanishing damping case, we need to make some modifications on Assumptions (A1)--(A3). The assumption corresponding to the divergence condition (\ref{eq:AssumAlpha}) is modified to
\begin{equation}\label{eq:divergence2}
\begin{split}
&\lim_{k\rightarrow \infty}\alpha_k=0,
\quad
\lim_{k\rightarrow \infty}\mu_k=0,
\quad
\lim_{k\rightarrow \infty}\frac{\alpha_k}{\mu_k}=0,
\quad\sum_{k=1}^\infty \alpha_k\mu_k=\infty,\\
& \exists L_{\mu} \geq 0,\quad \mu_{k-1}-\mu_{k}=L_{\mu} \alpha_{k} \mu_{k}+o\left(\alpha_{k} \mu_{k}\right).
\end{split}
\end{equation}
Let $\beta_k=\alpha_k/\mu_k$, then the $h_0$-slow condition corresponding to $\{\beta_k\}$ is that the divergence condition (\ref{eq:divergence2}) is satisfied and
\begin{equation}\label{eq:d-slower1}
 \exists K_s>0, \text{ such that }\quad \frac{\beta_n}{\beta_m} \ge K_s \prod_{k=m + 1}^n (1 - h_0 \alpha_k\mu_k) \quad \text{ for any }\ n > m \text{ large enough}.
\end{equation}
The sufficient decrease condition about $\{\beta_k\}$ is that the divergence condition (\ref{eq:divergence2}) is satisfied and
\begin{equation}\label{decrease2}
 \frac{\beta_k - \beta_{k+1}}{\alpha_{k}} = o(\alpha_{k}\mu_k)\quad \text{for large enough } k.
\end{equation}

By the asymptotic properties of (\ref{eq:divergence2}), for large enough $k$, the iteration (\ref{eq:mSGD4}) becomes
\begin{equation}\label{eq:mSGD5}
\begin{split}
Z_k = \left( I - \alpha_k \tilde{D} - \alpha_k\mu_k \tilde{E}+o(\alpha_k\mu_k)\right)Z_{k-1}  + \alpha_k
\left(\begin{array}{c}
0 \\
R_{k-1} + \xi_k
\end{array}\right)+o(\alpha_k\mu_k).
\end{split}
\end{equation}
Define $B_{n}=\prod_{k=1}^{n}\left(I-\alpha_{k} \tilde{D}-\alpha_{k} \mu_{k} \tilde{E}+o\left(\alpha_{k} \mu_{k}\right)\right)$.
 Timing $B_k^{-1}$ to both sides of \eqref{eq:mSGD5} and taking summation from $k=1$ to $n$,  we get
\begin{equation}
Z_n= B_n Z_{0}  + \sum_{k=1}^n B_n B_k^{-1} \left(\alpha_k
\left(\begin{array}{c}
0  \\
R_{k-1} + \xi_k
\end{array}\right)+o(\alpha_k\mu_k)\right).
\end{equation}

\begin{theorem}\label{thm:CLT-mSGD2} Suppose that  the divergence condition \eqref{eq:divergence2} and the $h_0$-slow condition~\eqref{eq:d-slower1}  hold with $h_0 < \mathrm{min}\{2/\zeta\mu,~2\}$, where $\zeta=1/(2L+2L_{\mu}^2)$, and there exists $K_\beta>0$, such that
\begin{equation}\label{eq:StepSizem1}
\beta_n:=\frac{\alpha_n}{\mu_n} \le K_\beta \min_{k \le n+1}\beta_k\qquad \text{ for all } n,
\end{equation}
and Assumptions (A6), (A9) hold, then we have
\begin{equation}\label{eq:CLT-Vanish}
\mathrm{Cov}(Z_n)^{-1/2}Z_n  \Rightarrow N(0,I).
\end{equation}
Furthermore, if $A \Sigma=\Sigma A$ and \eqref{decrease2} holds, there exists  a unique SPD matrix
 \begin{equation}
 W^{*}=\frac{1}{2}\left(\begin{array}{cc}A^{-1} \Sigma & 0 \\ 0 & \Sigma\end{array}\right)
 \end{equation}
 such that for mSGD \eqref{eq:mSGD4}, we have
\begin{equation}
\frac{Z_n}{\sqrt{\beta_n}} \Rightarrow N(0,W^*).
\end{equation}
\end{theorem}

The proof of Theorem \ref{thm:CLT-mSGD2} is similar to Theorem \ref{thm:CLT-vSGD}, relying on the following lemmas concerning the covariance and   $L^p$ bounds of $Z_n$,
and the proof of the following lemmas is similar
to that of mSGD. The proof of Lemma \ref{lm:mSGD1} is similar to Lemma \ref{lm:lm3} and we omit it.  The proof of Lemma \ref{lm:mSGD2} is deferred to Appendix \ref{sec:A5}.

 \begin{lemma}\label{lm:mSGD1}
Consider the triplet $(A_k, \beta_k, W)$ with the conditions
\[A_k  = \left(I-\alpha_{k} \tilde{D}-\alpha_{k} \mu_{k} \tilde{E}+o\left(\alpha_{k} \mu_{k}\right)\right) A_{k-1},\]
where  $\beta_k \rightarrow 0$, $W = W^T$.  Further assume  the learning rates satisfy condition~\eqref{eq:divergence2},  then we have
\begin{equation}
\sum_{k=1}^{n} A_n A_{k}^{-1}W A_k^{-T} A_n^T \alpha_k \mu_k\beta_k \rightarrow 0.
\end{equation}
\end{lemma}

\begin{lemma}\label{lm:mSGD2}
Assume that $f$ is $L$-smooth and $\mu$-strongly convex, and the Assumptions in Theorem~\ref{thm:CLT-mSGD2}  are made. We have the following upper bounds
\begin{equation}\label{eq:UBound1}
 \quad\mathbb{E}\| Z_n  \|^{2p} \le U_p \beta_n^p, \quad \text{ for }\quad p \in \Big[1, 1+\frac{\delta_0}{2}\Big],
\end{equation}
where $\delta_0$ is the constant defined in \eqref{eq:XiBound}.
\end{lemma}

\proof{Proof of Theorem~\ref{thm:CLT-mSGD2}}
Similar to Theorem \ref{thm:CLT-vSGD}, we first consider the linear case.
The mSGD with vanishing damping has the form
\begin{equation}\label{eq:z_nLinear}
Z_n= B_n Z_{0}  + \sum_{k=1}^n B_n B_k^{-1} (0,\alpha_k\xi_k^T)^T, ~ B_{n}=\prod_{k=1}^{n}\left(I-\alpha_{k} \tilde{D}-\alpha_{k} \mu_{k} \tilde{E}+o\left(\alpha_{k} \mu_{k}\right)\right).
\end{equation}
Since $B_nX_0\rightarrow 0$ almost surely, we will skip it in the analysis below.

Following similar steps in Theorem~\ref{thm:CLT-vSGD}, we only need slight modifications  in the estimation of covariance. Without loss of generality, assume there exists $\gamma> 0$, and $\Sigma_k\succ \gamma I$. Let $W=\operatorname{diag}(A^{-1}, I) / 2$,  then we have $W \tilde{D}+\tilde{D} W^{T}=0, W \tilde{E}+$ $\tilde{E} W^{T}=\tilde{E}$. From \eqref{eq:StepSizem1} and
\[
\begin{aligned}
B_{n} B_{k}^{-1} W B_{k}^{-T} B_{n}^{T}-B_{n} B_{k-1}^{-1} W B_{k-1}^{-T} B_{n}^{T} = B_{n} B_{k}^{-1} \tilde{E}  B_{k}^{-T} B_{n}^{T} \alpha_{k} \mu_{k} +  B_{n} B_{k}^{-1}\tilde{F}_k B_{k}^{-T} B_{n}^{T} \alpha_{k} \mu_{k},
\end{aligned}
\]
where $\tilde{F}_k = o(1)$, the summation of last term converges to $0$ by Lemma \ref{lm:mSGD1}.
Define $\tilde{\Sigma}_k=\left(\begin{array}{cc}
0 & 0\\
0 & \Sigma_k
\end{array}\right)$. Then we have
\[\begin{aligned}
\operatorname{Cov}\left(Z_{n}\right) & =  \sum_{k=1}^{n} B_{n} B_{k}^{-1} \tilde{\Sigma}_k B_{k}^{-T} B_{n}^{T}\alpha_{k}^{2}
 \succeq K^{-1}_{\beta} \gamma \sum_{k=1}^{n} B_{n} B_{k}^{-1} \tilde{E} B_{k}^{-T} B_{n}^{T} \alpha_{k} \mu_{k} \beta_n \\
& = K^{-1}_{\beta} \gamma \beta_n \sum_{k=1}^{n}\left( B_{n} B_{k}^{-1} W B_{k}^{-T} B_{n}^{T}-B_{n} B_{k-1}^{-1} W B_{k-1}^{-T} B_{n}^{T}\right)+o\left(\beta_n\right) \\
&=K^{-1}_{\beta} \gamma \beta_n\left(W-B_{n} W B_{n}^{T}\right)+o\left(\beta_n\right) \succeq C \beta_n,
\end{aligned}\]
which means $\operatorname{Cov}\left(Z_{k}\right) /\beta_{k}$ has lower bound.  The rest proof of \eqref{eq:CLT-Vanish} is similar to that of Theorem \ref{thm:CLT-vSGD}, which is omitted here.

To characterize the explicit form of the limit covariance matrix, we need commutativity between $A$ and $\Sigma$.
Let $V_{k}=\operatorname{Cov}\left(Z_{k}\right), W_{k}=V_{k} / \beta_{k+1}$, and $\Sigma=\lim\limits_{k\rightarrow\infty}\mathbb{E} \xi_{k} \xi_{k}^{T}$. Since $\beta_{k}-\beta_{k-1}=o\left(\alpha_{k} \mu_{k}\right)$, we have
\[
W_{k}-W_{k-1}=-\alpha_{k}\left(\tilde{D} W_{k-1}+W_{k-1} \tilde{D}^{T}\right)-\alpha_{k} \mu_{k}\left(\tilde{E} W_{k-1}+W_{k-1} \tilde{E}-\tilde{\Sigma}\right)+o\left(\alpha_{k} \mu_{k}\right),
\]
 which gives  $\tilde{D} W^{*}+W^{*} \tilde{D}^{T}=0$, $\tilde{E} W^{*}+W^{*} \tilde{E}+\tilde{\Sigma}=0$  if $W_{k}\rightarrow W^{*}$, where $\tilde{\Sigma}$ is defined as in Theorem \ref{thm:CLT-mSGD1}. Combing with the condition $A\Sigma=\Sigma A$, we get
 \[W^{*}=\frac{1}{2}\left(\begin{array}{cc}A^{-1} \Sigma & 0 \\ 0 & \Sigma\end{array}\right).\]
 Let $U_{k}=W_{k}-W^{*}$. We obtain
\[
U_{k}=\left(I-\alpha_{k} \tilde{D}+\alpha_{k} \mu_{k} \tilde{E}\right) U_{k-1}\left(I-\alpha_{k} \tilde{D}^{T}+\alpha_{k} \mu_{k} \tilde{E}^{T}\right)+o\left(\alpha_{k} \mu_{k}\right).
\]
Similar to the proof of Theorem~\ref{thm:CLT-vSGD}, we get $W_{k} \rightarrow W^{*}$.

For the nonlinear case, the analysis is analogous to Theorem \ref{thm:CLT-vSGD}.

\endproof

\begin{remark}
The conditions \eqref{eq:divergence2} and \eqref{eq:d-slower1} are reasonable
and $\left\{\alpha_{k}\right\}$, $\left\{\mu_{k}\right\}$ can be chosen such that they satisfy both conditions.
 For example, we can choose $\alpha_{k}=K / k^{a}, ~K>0, ~a \in(0,1), ~\mu_{k}=K_{\mu} / k^{b}, ~K_{\mu}>0, ~b \in(0, a), $ and $a+b\leq 1$. Furthermore, if we choose $\alpha_{k}=K / k^{a}, ~K>0, ~a \in(1 / 2,1)$ and $\mu_{k}=K_{\mu} / \sum_{k=1}^{n} \alpha_{k},~K_{\mu}>0$, this corresponds to $\mu(t) \sim O(1 / t)$ in the continuum limit, which is related to the Nesterov's SGD case \citep{NAG1}.
\end{remark}

{For mSGD with vanishing damping, it is different from two-time-scale stochastic approximation in \cite{konda2004convergence}.
Consider the two-time-scale stochastic approximation linear iterations of the form
\begin{equation}\label{eq:twoscale}
\begin{split}
w_k & = w_{k-1} - \alpha_k\mu_k (A_{11} w_{k-1} + A_{12}u_{k-1} + \xi_k), \\
u_{k} & = u_{k-1} - \alpha_k  (A_{21} w_{k-1} + A_{22}u_{k-1} + \eta_k),
\end{split}
\end{equation}
which has a two-time-scale asymptotic covariance
$$
\frac{\mathbb{E} w_kw_k^T}{\alpha_k\mu_k} \rightarrow \Sigma_{11},\ \frac{\mathbb{E} w_ku_k^T}{\alpha_k\mu_k} \rightarrow \Sigma_{12},\ \frac{\mathbb{E} u_ku_k^T}{\alpha_k} \rightarrow \Sigma_{22}.
$$
Now comparing with linear mSGD iteration (\ref{eq:mSGD2}) with vanishing damping $\mu_k \rightarrow 0$:
\begin{equation}\label{eq:mSGD5}
\begin{split}
x_{k}&= x_{k-1}+\alpha_{k} v_{k}, \\
v_{k}&= v_{k-1}-\mu_{k} \alpha_{k} v_{k-1}-\alpha_{k} A x_{k-1} +\alpha_{k} \xi_{k},
\end{split}
\end{equation}
we have $\mathbb{E} z_k z_k^T / \beta_k \rightarrow W^* > 0$.
The difference between the two-time-scale stochastic approximation and linear mSGD iteration is that the two-time-scale form does not care about fast variable $u_k$, while in mSGD with vanishing damping, $x$ and $v$ have a strong correlation.
In fact, mSGD can be combined with the two-time-scale stochastic approximation, which is a problem to be considered in future work.}

\section{CLT for the time average}\label{sec:average}

In this section, we consider the CLT for the time average of $x_k$. Such type of results are quite common in the Monte Carlo methods or the stochastic gradient Langevin dynamics \citep{Robert-MC,SGLD}. However, as we will show, the CLT for the time average does not hold in general for the current setup in SGD.

Assume $f$ is twice differentiable, $L$-smooth (\ref{eq:L-smooth}), and $\mu$-strongly convex (\ref{eq:muStrong}).
In order to investigate the CLT for the time average $\bar{x}_n$, we define
 \[\bar{X}_n:=\frac{\sum_{k=1}^n \alpha_k X_{k-1}} {\sum_{k=1}^n \alpha_k}=\bar{x}_n - x^*,\quad X_{k}:=x_{k}-x^* \]
and we can get
\begin{equation*}
\begin{split}
\sum_{k=1}^{n}\alpha_kX_{k-1} &  = \sum_{i=1}^{n} \alpha_k A^{-1}(\xi_k + R_k) + X_{0} - X_{n}
\end{split}
\end{equation*}
from the vSGD (\ref{eq:vSGDX}).

For the linear case with quadratic loss $f(x) = x^TAx/2$, we have $R_k=0$ and thus
\begin{equation}\label{avSGD1}
\begin{split}
 \sum_{k=1}^{n}\alpha_kX_{k-1}  =& \left(I - \prod_{i=1}^n(I - \alpha_i A)\right) A^{-1}X_{0} + \sum_{k=1}^{n} \alpha_k A^{-1} \xi_k  \\
 &  - \sum_{k=1}^{n}\prod_{j=k+1}^n (I - \alpha_j A)\alpha_k A^{-1} \xi_k = A^{-1}(P_{1n}+P_{2n}+P_{3n}),
\end{split}
\end{equation}
where $P_{1n}, P_{2n}, P_{3n}$ are corresponding terms appearing in the first equality after removing $A^{-1}$ in each term due to the commutativity between $A^{-1}$ and $I-\alpha A$.

Define
\[T_n=\sum_{k=1}^n \alpha_k,\quad S_n=\sum_{k=1}^n \alpha_k^2.\]
For the learning rates, we assume $\sum_{k=1}^\infty \alpha_k^2 = +\infty$ and Assumption \ref{Assum1:alpha}.  We note that otherwise the condition $\sum_{k=1}^\infty \alpha_k^2 < +\infty$ ensures $P_{2n}\rightarrow Z\sim O(1)$ in $L^2$,  upon assuming suitable bounds on $\xi_k$, thus no CLT holds in general. Below we will show that in linear case, the following rescaling of $\bar{X}_n$ is $O(1)$ and
\begin{equation}\label{eq:CLT-Average}
\frac{T_n}{\sqrt{S_n}}\bar{X}_n \Rightarrow N(0,A^{-1}\Sigma A^{-1}).
\end{equation}

First we note that $P_{1n},P_{3n}\sim o(\sqrt{S_n})$ since they are bounded in the limit $n\rightarrow \infty$. So only $P_{2n}$ matters. We can build the CLT \eqref{eq:CLT-Average} based on the Lindeberg condition for $Y_{nk} = \alpha_k \xi_k/\sqrt{S_n}$, $1\le k\le n$, which can be established by imposing the conditions like \eqref{eq:XiBound}: there exists $\delta_0 > 0$, such that $\sup_k \mathbb{E} \| \xi_k \|^{2+\delta_0}  < +\infty$. So  we get
\begin{equation}\label{eq:CLT-P2n}
\frac{P_{2n}}{\sqrt{S_{n}}} \Rightarrow N(0,\Sigma), ~\text{i.e.,}~ \frac{\sum_{k=1}^n \alpha_k A^{-1}\xi_k }{\sqrt{S_n}} \Rightarrow N(0,A^{-1}\Sigma A^{-1})
\end{equation}
and we are done.

However, the CLT for the time average does not hold for general nonlinear case, even when we consider the same rescaling of $\bar{X}_n$ as \eqref{eq:CLT-Average}. Below we will present a simple example to show this point.

Consider the loss function with strong convexity $f(x)=\mathbb{E}_\omega f(x,\omega)$, where
\[
f(x,\omega)=\frac{x^2}{2}+\phi(x)+b(\omega)x,~\phi(x)=\frac{x^3}{\sqrt{1+x^6}},~x\in\mathbb{R},
\]
and $b(\omega)$ is a bounded random variable with mean $0$. So we have $\xi_k=b(\omega_k)$ with $\mathbb{E}[ \xi_k|x_{k-1}] =0$ and it fulfills all necessary moment bounds condition.

The key point is to estimate the term $\sum_{k=1}^{n} \alpha_k R_k$. Take $\alpha_k=k^{-c}$, $c\in (0,1/2)$. At first, it is not difficult to show that $\mathbb{E}x_k^2\ge C\alpha_k$ by following similar approach as Step 2 in proving Theorem~\ref{thm:CLT-vSGD}, where the nonlinear terms can be controlled by the $L$-smoothness of $f$. Direct calculation gives $\phi^{\prime}(x)=3 x^{2}/(1+x^{6})^{\frac32}$, $\phi'''(0)\neq0$,  and $R(x)=f'(x)-x=O(x^2)$. We have
\[
R(x) \leq 3 x^{2} \text{ for any }x\in \mathbb{R}, \text{ and }R(x)\ge x^{2} \text{ for }|x|\le 1.
\]
 Similar to the estimation of the nonlinear term in Theorem \ref{thm:CLT-vSGD} with Lemma \ref{lm:vSGD}, we know that $\mathbb{E}R_k\leq 3\mathbb{E} x_{k-1}^{2} \leq C \alpha_{k}$
and $\mathbb{E} R_{k}\ge \mathbb{E} x_{k-1}^{2} 1_{\{|x_{k-1}|<1\}}-\mathbb{E} R_{k} 1_{\{|x_{k-1}|>1\}}$. Since $\mathbb{E} R_{k} 1_{\{|x_{k-1}|>1\}} \leq 3 \mathbb{E}|x_{k-1}|^{2+p} \leq C\alpha_{k}^{1+\frac{p}{2}}$ and $\mathbb{E}x_k^2\ge C\alpha_k$,
there exists $K>0$ such that  $\mathbb{E} R_{k} \geq K \alpha_{k}$. So we have
\begin{equation*}
\mathbb{E} \left(\frac{\sum_{k=1}^{n} \alpha_k R_k}{\sum_{k=1}^{n} \alpha_k^2}\right)  = \Theta(1),
\end{equation*}
which means that $T_n/\sqrt{S_n}\cdot\bar{X}_n\rightarrow\infty$, i.e., no CLT holds for the time average in this nonlinear case when we consider the same rescaling  as the linear case \eqref{eq:CLT-Average}.

Interestingly, the CLT  can be recovered for nonlinear case by assuming stronger conditions on $R(x)$.
We can show that under the condition $\sum_{k=1}^{+\infty}\alpha_k^2 = \infty$, Assumptions (A1) and (A4), and the learning rates $\{\alpha_k\}$ satisfy the $h_0$-slow condition, then there exists $C\ge 0$ such that
\[
\mathbb{E}\left(\frac{\sum_{k=1}^n \alpha_k X_{k-1}}{\sum_{k=1}^n \alpha_k^2}\right) \le C.
\]
Moreover, if  $ \| R(x)  \| \le C \|x\|^{2 + q},~ q > 0$, $Y_{nk} = \alpha_k \xi_k/\sqrt{S_n}$  fulfills the Lindeberg condition
and
\[
\frac{\sum_{k=1}^n \alpha_k^{2+q/2} }{\sqrt{\sum_{k=1}^n \alpha_k^2}} \rightarrow 0,
\]
 then we have the CLT \eqref{eq:CLT-Average}.
The analysis of this conclusion is almost the same as Theorem \ref{thm:CLT-vSGD}. For a specific example, we can choose a centrally symmetric function $f(x) \in C^4$ such that  $\| R(x)\| \le K \| x \|^{3}, q = 1$.  And for $\alpha_k = K n^{-a}$,  $a \in (1/4,1/2]$, the CLT exists.

\section{Numerical experiments}\label{sec:numer}

\begin{figure}[t]\label{fig:vSGD1}
\centering		
\begin{minipage}[c]{1\textwidth} 		
\centering			
\includegraphics[width=1\textwidth]{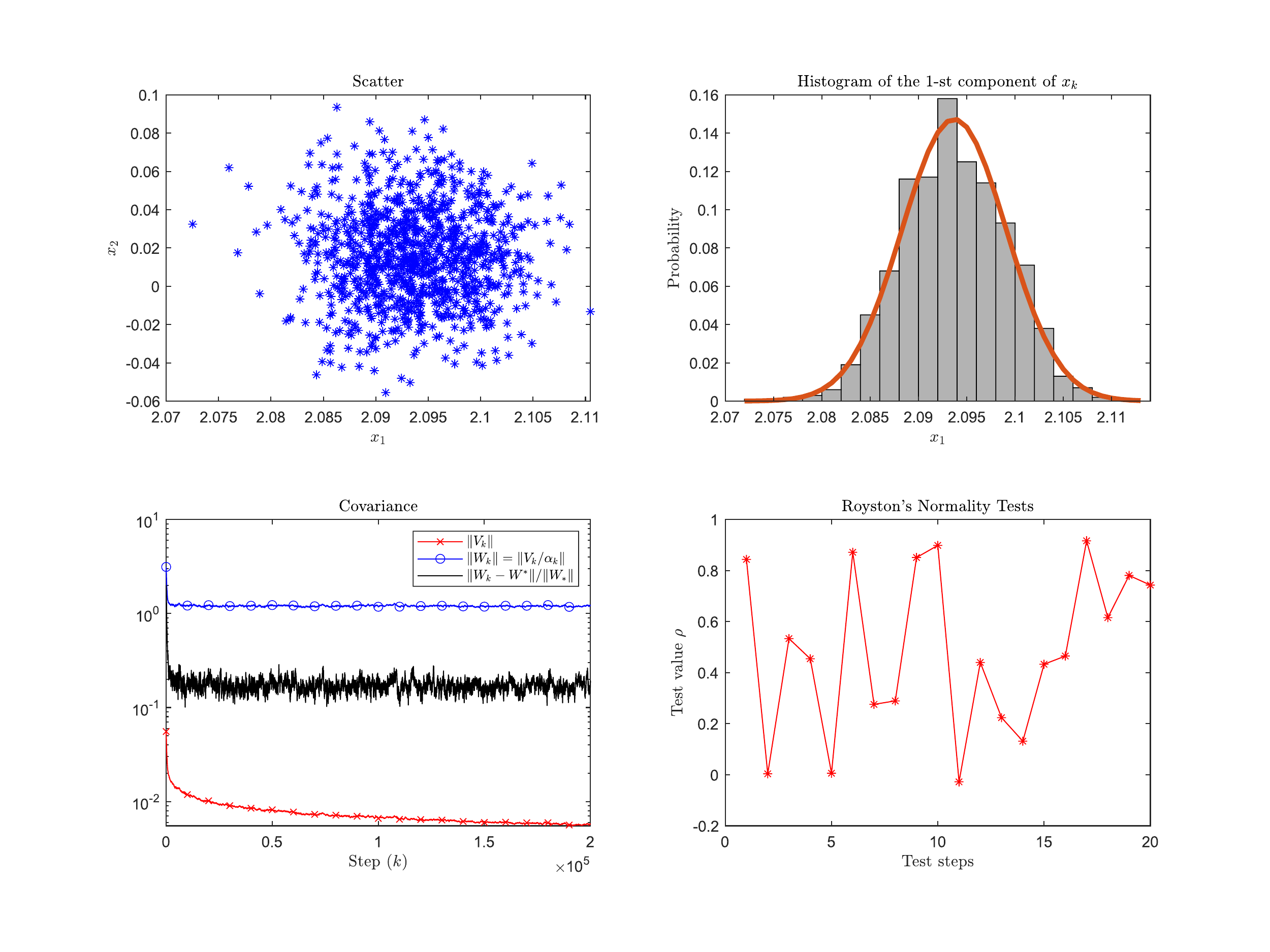}
\end{minipage}	
\caption{(Color Online). Numerical simulation of vSGD with $a = 0.25$.
Top left: Scatter plot of the first two components of $\{ x_k \}$ after $2\times 10^5$ steps.
Top right:
The distribution of the first component of $\{ x_k \}$ after $2\times 10^5$ steps.
Bottom Left:  $\|V_k\|=\| \mathbb{E}(x_k x_k^T)\|$ (red stars),  $\|W_k\|=\| \mathbb{E}(x_k x_k^T)/\alpha_k \|$ (blue circles) and relative error $\|W_k - W^*\|/\| W^* \|$ (black line).
Bottom right: Test values $\rho$ per $10^4$ steps. }
\label{figure:vSGD1}
\end{figure}

In this section, we will present numerical experiments for a toy model to validate the  theoretical analysis for the CLT  for the SGD-type methods.

We consider a logistic regression problem, i.e., the finite sum with penalty function, $\mu$-strongly convex function $$f(x)=  \frac{1}{N}\sum_{i=1}^n\left[ \ln(1 + e^{w_i^T x}) - y_i w_i^T x \right]+ \frac{\beta }{2}\| x\|^2, x \in \mathbb{R} ^d
$$
By computing, we have \[
\nabla f(x) = \frac{1}{N}\sum_{i=1}^n  \left[ \frac{w_i}{ 1 + e^{w_i^T x}} + (1 - y_i) w_i \right]+ \beta x, \nabla^2 f(x) = \frac{1}{N}\sum_{i=1}^n  \left[ \frac{w_i w_i^T e^{w_i^T x}}{ (1 + e^{w_i^T x})^2} \right] + \beta I.
\] 
In the following numerical experiments, we choose learning rates $\alpha_k =0.1 \cdot k^{-a}$ and $\beta=0.05$, $d = 10$, the batch size is 1.

\paragraph{Test of vSGD} We choose the sample size $M=1000$ in the simulation of vSGD with learning rates $\alpha_k = 0.1 \cdot k^{-0.25}$, and make statistical plots after $2\times 10^5$ steps (Figure \ref{figure:vSGD1}).

The top left panel of Figure \ref{figure:vSGD1} shows the scatter plot of $\{ x_k \}$ at the final simulation step.
The top right panel shows the distribution of the first component of $\{ x_k \}$ at the final simulation step, which perfectly fits a normal distribution. Other components of $\{ x_k \}$ shows similar behavior, which is omitted here. To check whether $W_k = \mathbb{E} x_k x_k^T / \alpha_k \rightarrow W^*$,  we first present the results to show how $W_k$ change with the number of iterations in the bottom left  panel of Figure \ref{figure:vSGD1}, which demonstrates convergence in a very early stage of iterations. To further quantify how far $W_k$ is from $W^*$, we utilize different sample sizes $M=100, 250, 500, 1000, 2000$ and learning schedules $\alpha_k=k^{-a}$ with $a=0.25, 0.50, 0.75$ to make the study, and the results are shown in Table  \ref{table:SGDtest}. To understand its meaning, we note that we indeed use $\tilde{W}_k=\sum_{m=1}^M x_{k}(\omega_m)x_{k}^T(\omega_m)/(M\alpha_k)$ to approximate $W_k = \mathbb{E} x_k x_k^T / \alpha_k$, and we have  the bias-variance error decomposition
\[\tilde{W}_k- W^*=\tilde{W}_k-W_k+W_k- W^*\sim  O(M^{-1/2}) + O(k^{a-1}).\]
When $a$ is relatively small ($a = 0.25,0.5$ in Table \ref{table:SGDtest}), the bias error can be neglected and we may expect the decay ratio $\sqrt{2}/2$ when $M$ is doubled. This can be exactly observed in the cases $a = 0.25,0.5$. However, when $a=0.75$, the bias error is not that small compared with the sampling error, and we can not get such perfect scaling behavior although the error decays with the increasing sample size.

\begin{figure}[ht]\centering		
\begin{minipage}[c]{1\textwidth} 		
\centering			
\includegraphics[width=1\textwidth]{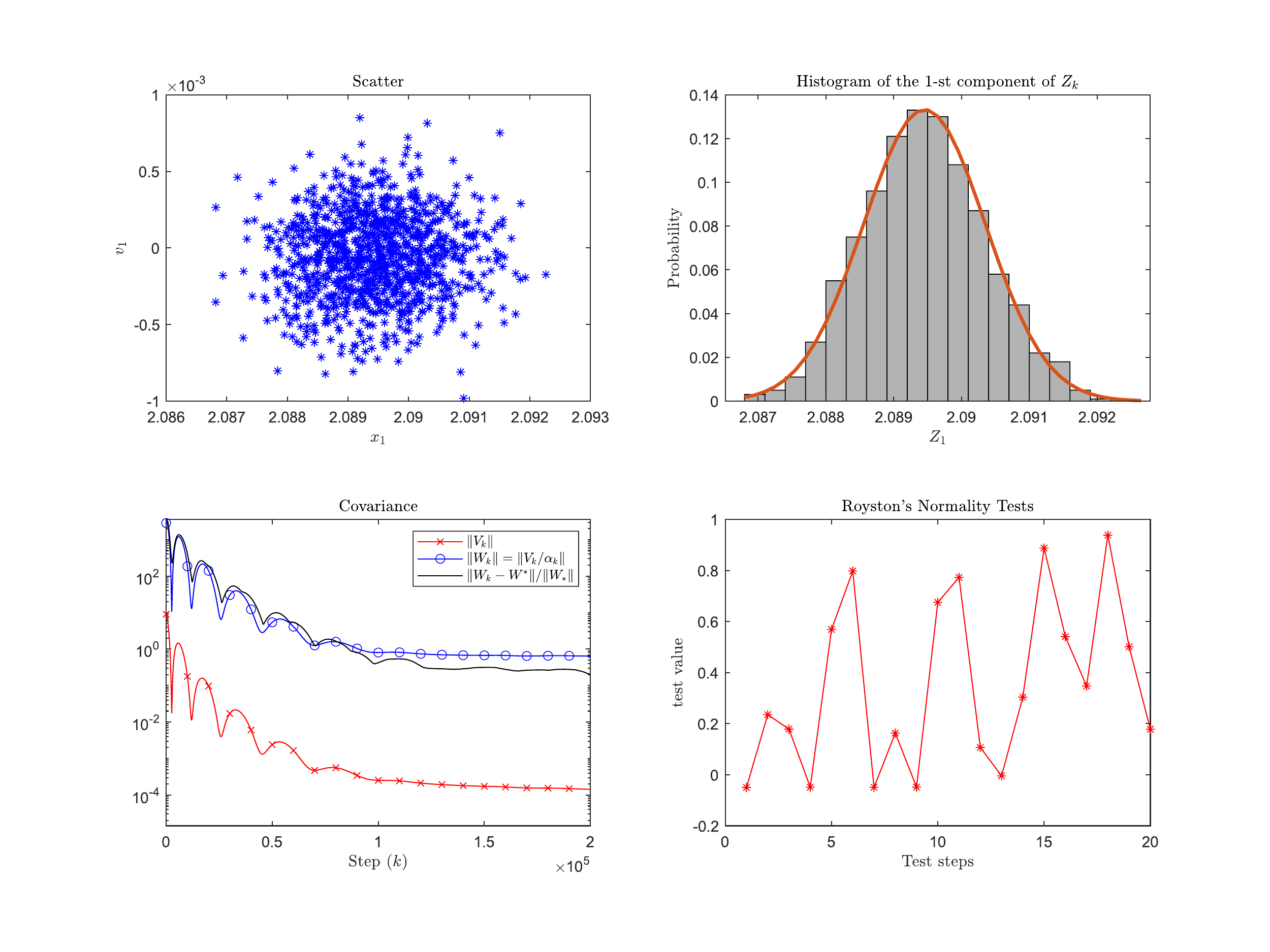}
\end{minipage}	
\caption{(Color Online). Numerical simulation of mSGD with $a = 0.5$ and $\mu = 0.2$.
Top left: Scatter plot of the first components of $\{ x_k \}$ and $\{ v_k \}$  after $2\times 10^5$ steps.
Top right: Distribution of the first  component of $\{ Z_k \}$ after $2\times 10^5$ steps.
 Bottom left:  $\|V_k\|=\| \mathbb{E}(Z_k Z_k^T)\|$ (red stars),  $\|W_k\|=\| \mathbb{E}(Z_k Z_k^T)/\alpha_k \|$ (blue circles) and relative error $\|W_k - W^*\|/\| W^* \|$ (black line). Bottom right: Test values $\rho$ per $10^4$ steps.
 }
\label{figure:mSGD0}
\end{figure}

\begin{table}[t]
\caption{Relative error $\|\tilde{W}_k-W^*\|/\|W^*\|$ for vSGD with different sample sizes $M$ and learning rates $\alpha_k=k^{-a}$.}
\centering
\begin{tabular}{| c|c|c|c|c|c|}
 \hline
  Sample size $M$ & 100 & 250 & 500 & 1000 & 2000 \\
  \hline
 $a = 0.25$ & 0.122 & 0.0761 & 0.0545 & 0.0385 & 0.0274 \\
 $a = 0.50$ & 0.123 & 0.0791 & 0.0535 & 0.0400 & 0.0268 \\
 $a = 0.75$ & 0.151 & 0.102 & 0.0957 & 0.0808 & 0.0646 \\
 \hline
\end{tabular}
\label{table:SGDtest}
\end{table}

To further quantitatively check $\mathrm{Cov}(x_k)^{-1/2}x_k \rightarrow N(0,1)$, we use Royston's multivariate normality test \citep{0Royston} to check whether $\{x_k\}$ is normally distributed at the $k$th simulation step. The test value is $\rho = p - \sigma$, where $p$ is the $p$-value associated with the Royston's statistic and $\sigma$ is the given significance ($\sigma=0.05$ in the experiments). The value $\rho>0$ indicates that $\{x_k\}$ passes the normal distribution test.
 The bottom right panel of Figure \ref{figure:vSGD1} shows the value of $\rho$ per  $10^4$ steps, which clearly verifies  $\mathrm{Cov}(x_k)^{-1/2}x_k \rightarrow N(0,1).$

\paragraph{Test of mSGD} We also choose $M=2000$ samples in the simulation of mSGD with learning  rates $\alpha_k = 0.1 \cdot k^{-a}$ to verify our CLT results after $5\times 10^5$ iteration steps (Figure~\ref{figure:mSGD0}).

We perform simulations using the mSGD method with constant damping $\mu = 0.2$ and learning rates $\alpha_k = 0.1\cdot k^{-0.5}$. The corresponding statistical plots are shown in Figure~\ref{figure:mSGD0}, and similar insights can be gained as the vSGD case.
The top left panel  shows the scatter plot of $\{ x_k \}$ and $\{ v_k \}$ at the final simulation step.
The top right panel  shows the distribution of the first component of $\{ Z_k \}$ at the final simulation step.
To verify $W_k = \mathbb{E} x_k x_k^T / \alpha_k \rightarrow W^*$,  we present the results to show how $W_k$ change with the number of iterations in bottom left panel of Figure \ref{figure:mSGD0}. Since the eigenvalues of $D$ are complex-valued, the oscillatory relaxation can be observed in the convergence process. The bottom right panel shows the value of $\rho$ per $10^4$ steps.  The positive values of $\rho$ suggest $\mathrm{Cov}(Z_k)^{-1/2}Z_k \rightarrow N(0,1).$

We also perform simulations using the mSGD method with vanishing damping $\mu_k = k^{-0.15}$ and learning rates $\alpha_k = 0.5\cdot k^{-0.75}$. The numerical results again support our theoretical analysis $\mathrm{Cov}(Z_k)^{-1/2}Z_k \rightarrow N(0,1)$, and they pass the normality test well. However, we will omit it due to the limit of space.

\section{Conclusion}

In this article, we re-established the CLT of vSGD under more general learning rates conditions. We also studied the mSGD and NaSGD with constant damping $\mu_k\equiv\tilde{\mu}$ and vanishing damping $\mu_k\rightarrow 0$ by taking advantage of the Lyapunov function technique, which are not previously known. The CLT for the time average was also investigated, and we found that  the CLT  for the time average held in the linear case, while  it was not true in general nonlinear situation. Numerical tests were carried out to verify the theoretical analysis. Applications and further extensions of the results obtained in this paper will be of future interest.

\section{Appendix: Proof of Lemmas in the main text}
\subsection{Proof of Lemma \ref{lm:lm3}}\label{sec:A1}

The proof of Lemma~\ref{lm:lm3} relies on the following lemma.
\begin{lemma} \label{lm:lmseq}
Assume the learning rates $\{\alpha_k\}$ satisfy Assumption (A1) and  the positive sequence $\beta_k\rightarrow 0$. Consider the triplet $(z_k, \beta_k, c )$ with the relation
\begin{equation}\label{eq:lemma41}
z_k \le (1 - c\alpha_k) z_{k-1} + \alpha_k \beta_k,\qquad z_k, c>0.
\end{equation}
We have: (a) $z_k \rightarrow 0$. Further assume there exists $M > 0$, such that $\beta_n \ge M\prod_{k=m + 1}^n (1 - h \alpha_k) \beta_m$,~for any $n > m$, where $\beta_k$ satisfies the $h$-slow condition with $h < c$. Then we have: (b) $z_n \le K \beta_n$, where  $ K>0$ is a constant.
\end{lemma}
\proof{Proof of Lemma \ref{lm:lmseq} }
To establish (a), we first assume $\sup_k \alpha_k < c^{-1}$ without loss of generality.
Let $B_n = \prod_{k=1}^n (1 - c\alpha_k)$,  we have
\[
z_n \le \sum_{k=1}^n B_n B_k^{-1} \alpha_k \beta_k+ B_n z_0.
\]
By $ B_n B_k^{-1} c\alpha_k  = B_n B_k^{-1} - B_n B_{k-1}^{-1}$,  we get $\sum_{k=1}^n  B_n B_k^{-1}\alpha_k \rightarrow c^{-1}.$

Let $\gamma_n = \sup_{k \ge n} {\beta_k}, $ then for large enough $m,n$ and $m < n$, we have
\[
 z_n  \le B_n B_m^{-1} z_m + \sum_{k=m+1}^n B_n B_k^{-1} \alpha_k \gamma_m \le c^{-1}\gamma_m +  B_n B_m^{-1} z_m.\]
 Taking $m = 0,$ we obtain that $\{z_k \}$ is bounded.
  Then  for any $ \varepsilon>0$, there exists $M>0$,  such that $\gamma_m < c\varepsilon/2$, for any $m > M$.
  And there exists $N>0$,  such that $B_n B_k^{-1} < \varepsilon /2$, for any $m > N$, which gives $z_k \rightarrow 0.$

For (b), letting $C_n = \prod_{k=1}^n (1 - c\alpha_k)(1 - h \alpha_k)^{-1}$, we get $\sum_{k=1}^n \alpha_k C_n C_k^{-1} \rightarrow (c-h)^{-1}$. By
 \[z_n \le \sum_{k=1}^n C_n  C_k^{-1} \alpha_k \beta_n+ B_n z_0\]
and $B_nz_0\rightarrow 0$, we obtain that there exists a constant $K>0$ such that $z_n \le K \beta_n$.
\endproof

The same result can be obtained if condition (\ref{eq:lemma41}) is replaced with
\begin{equation*}
z_k \le (1 - c\alpha_k + o(\alpha_k) ) z_{k-1} + \alpha_k \beta_k.
\end{equation*}

\proof{Proof of Lemma \ref{lm:lm3}}
 Consider the Lyapunov equation $D^{T} W+W D=Q$ with $Q\succ0, W\succ0$. There exists $c>0$, such that $ Q / 2\succ 2 c W$. We have
\[
\begin{aligned}
A_{m}^{-T} A_{n}^{T} W A_{n} A_{m}^{-1}
&=A_{m}^{-T} A_{n-1}^{T}
  \left(W-\alpha_{n}\left(D^{T} W+W D\right)+o\left(\alpha_{n}\right)\right)
 A_{n-1} A_{m}^{-1} \\
& \preceq A_{m}^{-T} A_{n-1}^{T}
   \left(W-\alpha_{n} Q+o\left(\alpha_{n}\right)\right) A_{n-1} A_{m}^{-1} \\
& \preceq A_{m}^{-T} A_{n-1}^{T} (W+o\left(\alpha_{n}\right)) A_{n-1} A_{m}^{-1}-A_{m}^{-T} A_{n-1}^{T} Q A_{n-1} A_{m}^{-1} \frac{\alpha_{n}}{2} \\
& \preceq \left(1-2 c \alpha_{n}\right) A_{m}^{-T} A_{n-1}^{T} W A_{n-1} A_{m}^{-1},
\end{aligned}
\]
where $m, n$ are large enough and $m<n$. We can get
\[
A_{m}^{-T} A_{n}^{T} W A_{n} A_{m}^{-1} \preceq \prod_{k=m+1}^{n}\left(1-2 c \alpha_{k}\right) W,
\]
which gives
\[
\begin{aligned}
\left\|A_{m}^{-T} A_{n}^{T} W A_{n} A_{m}^{-1}\right\|_2 \leq   K\prod_{k=m+1}^{n}\left(1-2 c \alpha_{k}\right),
\end{aligned}
\]
where $K=\lambda_{\max }(W)$. So we have
\[
\left\|\sum_{k=1}^{n} A_{n} A_{k}^{-1} W A_{k}^{-T} A_{n}^{T} \alpha_{k} \beta_{k}\right\|_2 \leq K \sum_{k=1}^{n} \prod_{i=k+1}^{n}\left(1-2 c \alpha_{i}\right) \alpha_{k}\left|\beta_{k}\right|.
\]
Let $z_{n}=\sum_{k=1}^{n} \prod_{i=k+1}^{n}\left(1-2 c \alpha_{i}\right) \alpha_{k}\left|\beta_{k}\right|$ and consider the triplet $\left(z_{k},\left|\beta_{k}\right|, 2 c\right)$ in Lemma \ref{lm:lmseq},  then we know that $z_{n} \rightarrow 0$, which means
\[
\sum_{k=1}^{\infty} A_{n} A_{k}^{-1} W A_{k}^{-T} A_{n}^{T} \alpha_{k} \beta_{k} \rightarrow 0.
\]
The proof is completed.
\endproof

\subsection{Proof of Lemma \ref{lm:vSGD}}\label{sec:A2}

\proof{Proof of Lemma \ref{lm:vSGD}}
We will separate the cases $p=1$ and $p>1$.

\textit{Case $p=1$}.  By the $L$-smoothness and $\mu$-strong convexity of $f$,  we have
\begin{equation*}
\begin{split}
\mathbb{E} \| X_k \|^2 &\le \mathbb{E} \| X_{k-1} \|^2 - 2\alpha_k X_k^T\nabla f(x_{k-1}) + \alpha_k^2 (1+K_\xi) \| \nabla f(x_{k-1}) \|^ 2 + M\alpha_k^2 \\
&\le \big(1 - 2\mu \alpha_k + (1+K_\xi)L^2\alpha_k^2 \big)\mathbb{E} \| X_{k-1} \|^2  +  M\alpha_k^2.
\end{split}
\end{equation*}
Now take the triplet $( \mathbb{E} \| X_k \|^2, M\alpha_k, 2\mu )$ in Lemma \ref{lm:lmseq}. We get  $\mathbb{E} \| X_k  \|^2 \rightarrow 0.$  And with the $h_0$-slow condition of  $\alpha_k$,  we know that there exists $C>0$, such that $\mathbb{E} \| X_k  \|^2 \le C \alpha_k$.

\textit{Case $p>1$}. Define $Y_{k-1} = X_{k-1} - \alpha_k \nabla f(x_{k-1})$.  Through $\mu$-strongly convex and $L$-smooth condition of $f$, we have
\begin{equation*}
\begin{split}
\| Y_{k-1}\|^{2p} & = \big(\|X_{k-1}\|^2 - 2\alpha_k X_{k-1}^T\nabla f(x_{k-1}) + \alpha_k ^2\|\nabla f(x_{k-1})\|^2\big)^{p}  \\
& \le \|X_{k-1}\|^{2p}(1  - 2\mu \alpha_k  +  L^2 \alpha_k^2 )^p.
\end{split}
\end{equation*}
For $p\le 1+\frac{\delta_0}{2}$, we have
 \[\mathbb{E}\|\xi_k \|^{2p} \le M_p + V_p  \mathbb{E}\| \nabla f(x_{k-1}) \|^{2p}\le M_p + V_p L^{2p} \mathbb{E}\| X_{k-1} \|^{2p}\]
 by Assumption \eqref{eq:XiBound}. By the inequality $(a+b)^p\le a^p + 2^{p-1}p(a^{p-1}b + b^{p})$ for $p>1$ and $a,b>0$, we get
\begin{equation*}
\begin{split}
\quad \mathbb{E}\| X_k \|^{2p}&  =  \mathbb{E} \| Y_{k-1} + \alpha_k \xi_k \|^{2p} =  \mathbb{E} | Y_{k-1}^T Y_{k-1} +2\alpha_k \xi_k^T Y_{k-1} + \alpha_k^2 \xi_k^T  \xi_k |^{p}\\
& \le   C\big\{ \mathbb{E}(\|\xi_k \|^2 | Y_{k-1}^T Y_{k-1} +2\alpha_k \xi_k^T Y_{k-1} |^{p-1})\alpha_k^{2} + \mathbb{E}\| \xi_k\|^{2p}\alpha_k^{2p}\big\}\\
&~~~~ +  \mathbb{E} | Y_{k-1}^T Y_{k-1}+2\alpha_k \xi_k^T Y_{k-1}|^p.
\end{split}
\end{equation*}
Define $q = p/(p-1)$. By using H\"{o}lder's inequality, for $\varepsilon_1 > 0$,  we have
\begin{equation*}
\mathbb{E} \big(\|\xi_k \|^2  | Y_{k-1}^T Y_{k-1} +2\alpha_k \xi_k^T Y_{k-1}|^{p-1}\big)\alpha_k^{2} \le \alpha_k\Big(\frac{\varepsilon_1}{q}\mathbb{E} | Y_{k-1}^T Y_{k-1}  +2\alpha_k \xi_k^T Y_{k-1}|^{p}+\frac{\alpha_k^{p}}{p\varepsilon_1^{p-1}} \mathbb{E} \| \xi_k\|^{2p}\Big),
\end{equation*}
which means
\begin{equation*}
\mathbb{E}\| X_k \|^{2p} \le (1 + C\varepsilon_1\alpha_k)\mathbb{E} | Y_{k-1}^T Y_{k-1} +2\alpha_k \xi_k^T Y_{k-1}|^{p} + C\alpha_k^{p+1}\mathbb{E}\| \xi_k\|^{2p}.
\end{equation*}
By H\"{o}lder's inequality and $(a+b)^p\le a^p+pa^{p-1}b + D_p(a^{p-2}b^2+b^p)$ for $p>1$ and $a,b>0$, where $D_p=2^{p-2}p(p-1)$,  for any $\varepsilon_2,~\varepsilon_3 > 0,$ we have
\begin{equation*}
\begin{split}
&\quad \ \mathbb{E} \big| Y_{k-1}^T Y_{k-1} +2\alpha_k \xi_k^T Y_{k-1}\big|^{p}  \\
& \le  \mathbb{E}\| Y_{k-1} \|^{2p} + p \mathbb{E} (\xi_k^T Y_{k-1} \| Y_{k-1} \|^{2p-2}) + D_p\big(\alpha_k^2 \mathbb{E} (\|\xi_k\|^{2}\| Y_{k-1} \|^{2p-2}) + \alpha_k^p \mathbb{E}(\| \xi_k\|^{p}\| Y_{k-1} \|^{p})\big) \\
& \le \mathbb{E} \| Y_{k-1} \|^{2p} + C\alpha_k \left(\mathbb{E}\frac{\alpha_k^p \| \xi_k\|^{2p}}{\varepsilon_2^{p-1}}  + \varepsilon_2\mathbb{E}\| Y_{k-1} \|^{2p}\right) + C \alpha_k^{p-1}  \left(\mathbb{E}\frac{\alpha_k^p \| \xi_k\|^{2p}}{\varepsilon_3}  + \varepsilon_3\mathbb{E}\| Y_{k-1} \|^{2p}\right)   \\
& \le  (1 +  C(\varepsilon_2 + \varepsilon_3)\alpha_k) \mathbb{E} \| Y_{k-1} \|^{2p} + C\alpha_k^{p+1}\mathbb{E}\| \xi_k\|^{2p}.
\end{split}
\end{equation*}
Then we get
\begin{equation*}
%\begin{split}
\mathbb{E}\| X_k \|^{2p} \le  (1 + C\varepsilon_1\alpha_k) (1 +  C(\varepsilon_2 + \varepsilon_3)\alpha_k) (1  - 2\mu \alpha_k + L^2 \alpha_k^2 )^p \mathbb{E}\| X_{k-1} \|^{2p} + C(1 +\mathbb{E}\| X_{k-1} \|^{2p})\alpha_k^{p+1}.
%\end{split}
\end{equation*}
By the arbitrariness of $\varepsilon_1, \varepsilon_2, \varepsilon_3$, we have
\begin{equation*}
\begin{split}
\mathbb{E}\| X_k \|^{2p} \le (1  - (2\mu - \varepsilon) p  \alpha_k + o(\alpha_k)) \mathbb{E}\| X_{k-1} \|^{2p} + C\alpha_k^{p+1}
\end{split}
\end{equation*}
for any $0<\varepsilon<2\mu$.
Now take the triplet $(\mathbb{E}\| X_k \|^{2p}, C\alpha_k^p, 2\mu - \varepsilon )$ in Lemma \ref{lm:lmseq}, we have
\begin{equation*}
\begin{split}
\mathbb{E}\|X_k  \|^{2p} \le U_p \alpha_k^p
\end{split}
\end{equation*}
with $U_p > 0$,  which gives the upper bounds of $X_k$ in $L^{2p}$ for $p\in [1,1+\delta_0/2]$.
\endproof

\subsection{$L^{2p}$ bounds for mSGD with constant damping $\mu_k\equiv\tilde{\mu}$} \label{sec:A3}

\proof{Proof}
For mSGD, we first consider the case $p=1$. Define the Hamiltonian $H_k = f(x_k)-f(x^*) + \|v_k\|^2/2 $ and  $\bar{H}_k = \| \nabla f(x_k)\|^2 + \|v_k\|^2$.

By \eqref{eq:mSGD31} and $L$-smoothness of $f$  and $\mathbb{E}[  \| \xi_k \|^2 |\mathcal{F}_{k-1}] \le M + K_{\xi} \| \nabla f(x_{k-1})\|^2$, we have
\begin{equation*}
\begin{split}
\mathbb{E} [f(x_k)|\mathcal{F}_{k-1} ] & \le f(x_{k-1}) + \alpha_k \mathbb{E} [\nabla f(x_{k-1})^T v_k|\mathcal{F}_{k-1} ] + \frac{L\alpha_k^2}{2}\mathbb{E} [ \| v_k\|^2|\mathcal{F}_{k-1} ]   \\
& \le f(x_{k-1}) + \alpha_k\nabla f(x_{k-1})^Tv_{k-1} + C\alpha_k^2 (1 + O (\bar{H}_{k-1})),
\end{split}
\end{equation*}
and
\begin{equation*}
\begin{split}
\mathbb{E}\Big[\frac{\| v_k \|^2}{2}|\mathcal{F}_{k-1}\Big] & \le \frac{\| v_{k-1}\|^2}{2} + \frac{1}{2}\mathbb{E} [ \|v_k - v_{k-1} \|^2|\mathcal{F}_{k-1} ] - \alpha_k \left(\tilde {\mu}\|v_{k-1}\|^2 + \nabla f(x_{k-1})v_{k-1}\right) \\
& \le \frac{\| v_{k-1} \|^2}{2} - \alpha_k (\tilde {\mu}\|v_{k-1}\|^2 + \nabla f(x_{k-1})v_{k-1}) + C\alpha_k^2(1 + O (\bar{H}_{k-1})),
\end{split}
\end{equation*}
thus
\begin{equation*}
\begin{split}
\mathbb{E} [H_k|\mathcal{F}_{k-1} ] & \le H_{k-1} -  \tilde{\mu} \alpha_k \| v_{k-1} \|^2 + C\alpha_k^2 (1 +O (\bar{H}_{k-1})).
\end{split}
\end{equation*}

By introducing the term $\tilde{Z}_k = v_k^T\nabla f(x_k)$, with $L$-smoothness of $f$, we get
\begin{equation*}
\begin{split}
\mathbb{E} [\tilde{Z}_k|\mathcal{F}_{k-1} ]  & = \tilde{Z}_{k-1} + \mathbb{E} [ v_{k}^T(\nabla f(x_k) - \nabla f(x_{k-1}))|\mathcal{F}_{k-1} ]+ \mathbb{E}[(v_k - v_{k-1})^T\nabla f(x_{k-1})|\mathcal{F}_{k-1} ] \\
&\le \tilde{Z}_{k-1} + \alpha_k (L \mathbb{E}[ \| v_{k} \|^2 |\mathcal{F}_{k-1} ]  -  \| \nabla f(x_{k-1}) \|^2  - \tilde{\mu} v_{k-1}^T \nabla f(x_{k-1})) \\
& = \tilde{Z}_{k-1} + \alpha_k (L \| v_{k-1} \|^2  -  \| \nabla f(x_{k-1}) \|^2  - \tilde{\mu} v_{k-1}^T \nabla f(x_{k-1})) \\
&\quad+ \alpha_k^2(O(\bar{H}_{k-1})+o(1)) .
\end{split}
\end{equation*}

Now we consider the Lyapunov function $H^E_k = \mathbb{E}\tilde{H}_k$ with $\tilde{H}_k = H_k + \zeta \tilde{Z}_k$, where $\zeta > 0$ is small enough.
From the $L$-smoothness and strong convexity of $f$, we have $\tilde{H}_k=\Theta (\bar{H}_k)$. Then
\begin{equation*}
\begin{split}
&\mathbb{E} [\tilde{H}_k|\mathcal{F}_{k-1} ]  \le \tilde{H}_{k-1} - \alpha_k F_{k-1} + C \alpha_k^2 (1 + O (\tilde{H}_{k-1})),
\end{split}
\end{equation*}
where $F_{k-1}  =  \zeta\| \nabla f(x_{k-1}) + \frac{\tilde{\mu}}{2}v_{k-1} \|^2 + (\tilde\mu - \zeta(\frac{\tilde{\mu}^2}{4} + L)) \| v_{k-1} \|^2$.
Since  $f$ is strongly convex, it is easy to get $F_k = \Theta(\bar{H}_k)$, which means that there exists $K>0$ such that
\begin{equation*}
\begin{split}
\mathbb{E} [\tilde{H}_k|\mathcal{F}_{k-1} ]  &\le (1 - K \alpha_k + o(\alpha_k) )\tilde{H}_{k-1} + C \alpha_k^2.
\end{split}
\end{equation*}
Then we have
\begin{equation*}
\begin{split}
H^E_k  &\le (1 - K \alpha_k + o(\alpha_k) )H^E_{k-1} + C \alpha_k^2.
\end{split}
\end{equation*}
Taking the triplet as  $( H^E_k, C\alpha_k, K )$ in  Lemma \ref{lm:lmseq}, we have $H^E_k \le U_1 \alpha_k$,
where $ U_1 > 0$. Then  the $L^2$ bound for mSGD is obtained as a result of $\mathbb{E}\| Z_k \|^2  \le C H^E_k$.

The estimate for $p \in [1,1+\delta_0/2]$ is similar to the derivations in Appendix~\ref{sec:A2} for vSGD and we omit it.
\endproof

\subsection{$L^{2p}$ bounds for NaSGD with constant damping $\mu_k\equiv\tilde{\mu}$}\label{sec:A4}

\proof{Proof} The analysis for NaSGD is similar to mSGD. For $p=1$, we consider the Hamiltonian
$H_k = f(x_k)-f(x^*) + \|v_k\|^2/2 $. Let $\bar{H}_k = \| \nabla f(x_k)\|^2 + \|v_k\|^2$. By (\ref{eq:NaSGD3}) with $y_{k-1} = x_{k-1} +  \beta_k v_{k-1}$, we have
\begin{equation*}
\begin{split}
\mathbb{E} [H_k|\mathcal{F}_{k-1} ] & \le H_{k-1} +\alpha_k [v_{k-1}^T(\nabla f(x_{k-1}) - \nabla f(y_{k-1})) - \tilde{\mu} \| v_{k-1}\|^2 ] + o(\alpha_k ) \\
& \le H_{k-1} -  (\tilde{\mu} + \mu\beta_k) \alpha_k  \| v_{k-1}\|^2 + C\alpha_k^2(1 +O (\bar{H}_{k-1})).
\end{split}
\end{equation*}
And for $\tilde{Z}_k =  v_k^T\nabla f(x_k) $, we have
\begin{equation*}
\begin{split}
\mathbb{E} [\tilde{Z}_k|\mathcal{F}_{k-1} ] &\le Z_{k-1} +   \alpha_k (L\| v_{k-1}\|^2 -   \nabla f(x_{k-1})^T\nabla f(y_{k-1}) - \tilde{\mu} v_{k-1}^T \nabla f(x_{k-1}) )+ \alpha_k^2 O (\bar{H}_{k-1}).
\end{split}
\end{equation*}
By $L$-smoothness of $f$, we get
\begin{equation*}
|\nabla f(x_{k-1})^T(\nabla f(x_{k-1})- \nabla f(y_{k-1}))| \le L\beta_k | v_{k-1}^T \nabla f(x_{k-1})| \le \hat{\beta}L \lambda \| \nabla f(x_{k-1}) \|^2 + \frac{\hat{\beta}L}{\lambda} v_{k-1}^2,
\end{equation*}
 where $\hat{\beta} = \sup_k \beta_k$ and $\lambda= 1/(2L\hat{\beta})$. Then we have
\begin{equation*}
\begin{split}
\mathbb{E} [\tilde{Z}_k|\mathcal{F}_{k-1} ] &\le Z_{k-1} +  \alpha_k \Big((L+ \frac{L\hat{\beta}}{\lambda})\| v_{k-1}\|^2 -  (1 - L\hat{\beta}\lambda) \| \nabla f(x_{k-1}) \|^2 - \tilde{\mu} v_{k-1}^T \nabla f(x_{k-1})\Big)\\
&+ C\alpha_k^2(1 + O (\bar{H}_{k-1})).
\end{split}
\end{equation*}
Now we consider the Lyapunov function $H^E_k = \mathbb{E}\tilde{H}_k=\mathbb{E}(H_k + \zeta \tilde{Z}_k)$ with small enough $\zeta>0$. Similar to analysis of mSGD in Appendix~\ref{sec:A3}, we have
\[
H^E_k \le (1 - K\alpha_k + o(\alpha_k))H^E_{k-1} + C\alpha_k^2.
\]
The rest analysis is similar to  the derivations in Appendix~\ref{sec:A3}.
\endproof

\subsection{ Proof of Lemma~\ref{lm:mSGD2}}\label{sec:A5}

\proof{Proof of Lemma~\ref{lm:mSGD2}}
For mSGD with $\mu_k \rightarrow 0$, consider the Lyapunov function $H^E_k = \mathbb{E}\tilde{H}_k$, where $\tilde{H}_k = H_k + \lambda \mu_k \tilde{Z}_k,$ with $0<\lambda < 1/L$. From the $L$-smoothness and strong convexity of $f$, we have $\tilde{H}_k=\Theta (\bar{H}_k)$. Similar to mSGD, we have
\begin{equation*}
\begin{split}
\mathbb{E} [\tilde{H}_k|\mathcal{F}_{k-1} ]  &\le \tilde{H}_{k-1} - \alpha_k \tilde{F}_{k-1} - \lambda(\mu_k - \mu_{k-1}) v_{k-1}^T \nabla f(x_{k-1})  +  C(1+O (\bar{H}_{k-1}))\alpha_k^2,
\end{split}
\end{equation*}
where $
\tilde{F}_{k-1} = \lambda \mu_k (\| \nabla f(x_{k-1}) \|^2 + v_{k-1}^T \nabla f(x_{k-1})) + (1 - L \lambda)\mu_k \| v_{k-1} \|^2.$

Further, let $\lambda < (L + L_\mu^2 / 4)^{-1}$. By (\ref{eq:divergence2}), we have $\tilde{F}_k + \lambda(\mu_{k+1} - \mu_{k})\tilde{Z}_{k} = \mu_{k+1} \Theta(\bar{H}_k)$.  Similar to the analysis of mSGD in Appendix~\ref{sec:A3}, there exists $K>0$ such that
\begin{equation*}
\begin{split}
\mathbb{E} [\tilde{H}_k|\mathcal{F}_{k-1} ]  &\le (1 - K \alpha_k \mu_k + o(\alpha_k \mu_k))\tilde{H}_{k-1} + C(1+O (\tilde{H}_{k-1}))\alpha_k^2. \\
\end{split}
\end{equation*}
Then we have
 \[H^E_k \le (1 - K \alpha_k\mu_k + o(\alpha_k \mu_k))H^E_{k-1} + C\alpha_k^2.\]
 Taking the triplet as  $(H^E_k, C\alpha_k / \mu_k, K)$ and using $\{ \alpha_k \mu_k\}$ instead of $\{\alpha_k \}$ in  Lemma \ref{lm:lmseq}, we have $H^E_k \le C \alpha_k / \mu_k$.

The estimate for $p \in [1,1+\delta_0/2]$ is also similar to the derivations in Appendix~\ref{sec:A2}.
\endproof

\section*{Acknowledgments}
The authors are grateful to the anonymous referees for careful reading and helpful suggestions.

\bibliographystyle{siam}
\bibliography{CLT-Ref}  

\end{document}